\documentclass{article}
\pdfoutput=1
\usepackage{amsfonts,amsmath}
\usepackage{makeidx}
\usepackage{graphicx}
\usepackage{multicol}
\usepackage{color}
\usepackage{listings}
\usepackage{geometry}

\title{Optimization Under Uncertainty Using the Generalized Inverse Distribution Function}

\author{
Domenico Quagliarella\\
Department of Fluid Mechanics\\
Italian Aerospace Research Center\\
Via Maiorise snc, 81043 Capua, Italy\\
d.quagliarella@cira.it
\and
Giovanni Petrone\\
Aerospace, Automotive and Turbo CFD Team\\ 
ANSYS UK Ltd,
Sheffield Business Park\\
6 Europa View, Sheffield, S9 1XH, UK\\ 
giovanni.petrone@ansys.com
\and
Gianluca Iaccarino\\
Mechanical Engineering Department\\
Stanford University, Stanford, CA, 94305, USA\\
jops@stanford.edu
}

\begin{document}

\maketitle

\abstract{A framework for robust optimization under uncertainty based
on the use of the generalized inverse  distribution function (GIDF),
also called quantile function, is here proposed. Compared to more classical
approaches that rely on the usage of statistical moments as deterministic
attributes that define the objectives of the optimization process, the
inverse cumulative distribution function allows for the use of all the
possible information available in the probabilistic domain. Furthermore,
the use of a quantile based approach leads naturally to a multi-objective
methodology which allows an a-posteriori selection of the candidate
design based on risk/opportunity criteria defined by the designer.
Finally, the error on the estimation of the objectives due to the
resolution of the GIDF will be proven to be quantifiable}

\section{Introduction}

Numerical design optimization procedures commonly imply that all the
design parameters can be precisely determined and that the manufacturing
process is reliable and exactly and indefinitely replicable, so that
it produces identical structures. Furthermore, it is not made any
assumption on the reliability and fidelity of the physical model that governs
the behavior of the product that is being designed.  Unfortunately,
industrial manufacturing processes and real operating conditions introduce
tolerances in the product and uncertainties in the working conditions that
may produce significant deviations from the conditions considered in the
design stage.  Robust optimization techniques are designed, conversely,
to try to overcome these problems and to account for uncertainty sources
since the numerical design optimization stage to avoid that discrepancies
between calculated and real performances may lead to a product that,
in the end, is not suitable for the purpose for which it was designed.

Another important source of discrepancy between the design and the
manufactured item originates from the fidelity of the physical model used
in the numerical design process.  This lack of physical knowledge is not
classifiable, in a strict sense, as a source of uncertainty, because,
rather than being statistic, it is epistemic, i.e. implicit in the
nature of the  computational model used.  Nevertheless, a credibility
analysis of deterministic results may be very helpful to improve the
numerical design process and the introduction of techniques inherited
from uncertainty quantification and robust design to take into account
these computational model inaccuracies, may be very helpful to improve
the quality and robustness of the resulting design.

In this work an approach to robust design optimization based on the use
of the generalized inverse distribution function is presented.  The robust
optimization framework is illustrated and the commonly used techniques to
face the problem are briefly summarized making reference to the related
literature.  The new approach is then introduced and illustrated with
the help of some examples built on top of mathematical test functions.
A very simple evolutionary multi-objective optimization algorithm based on
the usage of the inverse cumulative distribution function is illustrated
and, finally, some conclusive notes and remarks are drawn.

This work shares the same philosophy of another robust optimization
approach presented by the authors and based on the use of the cumulative
distribution function together with a template function that represents
a target ideal CDF. The CDF and the template are used to define a
``robustness index'' (RI) that measures the deviation from the given
optimal distribution \cite{Petrone2011b}. In the present work, the main
difference, as it will be possible to appreciate in the following, is
that the use of the inverse distribution function (also called quantile
function) allows to avoid the introduction of the robustness index.

\section{Robust optimization}

Let $Z$ be a metric space and $z \in Z$ the vector of design variables.
Let also $X:\Omega \rightarrow \Xi \subseteq \mathbb{R}$ be a real valued
random variable defined in a given $(\Omega,{\cal F}, P)$ probability
space.  We want to deal with an optimization problem where an objective
is optimized with respect to $z \in Z$ and depends on the realizations
$x$ of $X$.  In other terms we have:
$$
y(z, X): z \in Z, X 		\longrightarrow  Y(z)
$$
with $Y(z)$ a new random variable, e.g. a new mapping of $(\Omega,{\cal
F}, P)$ into $\mathbb{R}$, that depends on $z$.  Solving an optimization
problem involving $Y(z)=y(z,X)$ means that we want to find a value
$\bar{z} \in Z$ such that the random variable $Y(\bar{z})$ is optimal.
To establish the optimality of a given $Y(\bar{z})$ with respect to all
$Y(z),\; \forall z \in Z$, a ranking criterion must be defined such that
for any couple $z_1, z_2 \in Z$ it is possible to state that $Y(z_1)$
is better or worse than $Y(z_2)$ (from now on, $Y(z_1) \preceq Y(z_2)$
will mean that $Y(z_1)$ is better or equivalent to $Y(z_2)$).

Recalling that a random variable is a measurable function, it seems
natural to introduce measures that highlight particular features of the
function. This leads to the classical and widely used approach of using
the statistical moments to define the characteristics of the probability
distribution that are to be optimized.  More generally, let's consider
an operator
$$
\Phi_{X} : 	Y(z)=y(z, X) \in Z\times(\Omega,{\cal F}, P) \longrightarrow \Phi(z) \in V \subseteq \mathbb{ R}
$$
that translates the functional dependency on the random variable, $Y$,
into a real valued function of $z$ that represents a deterministic
attribute of the function, $Y(z)$.  This makes possible to formulate
the following optimization problem
$$
P_\Phi : \quad \min_{z \in Z} \Phi(z)
$$
Without loss of generality, it is possible to identify the random
variable $Y$ through its distribution function $f_Y(y)$ or its cumulative
distribution function $F_Y(y)$. If $\Phi(\cdot)$ is assumed as the expected
value of the objective function ($\mathbb{E}$), the classical formulation
of first moment optimization is retrieved:
$$
P_{_\mathbb{E}} : \quad \min_{z \in Z} \left.\int_{\mathbb{ R}} y f_Y(y,z)dy\right.
$$
that in terms of the CDF becomes:
$$
P_{_\mathbb{E}} : \quad \min_{z \in Z} \left.\int_{\mathbb{ R}} y d F_Y(y,z)\right.
$$
It should be noted that here the distribution function depends also on
$z$, that is the vector of the design variables.

For the purposes of the definition of the problem, it is not necessary to
know exactly the distribution $f_Y$ (or $F_Y$). Indeed, it is possible,
as will be seen below, to use an estimate of the distribution having the
required accuracy. In particular, the Empirical Cumulative Distribution
Function (ECDF) will be used in this work as statistical estimator of
the CDF.

The first order moment method is also called mean value approach, as
the mean is used as objective to reduce the dependency on $Y$. This
method is widely used, mostly because the mean is the faster converging
moment and relatively few samples are required to obtain a good
estimate. Often, however, the mean alone is not able to capture
and represent satisfactorily the uncertainties embedded in a given
design optimization problem.  Tho overcome this drawback, a possible
approach is the introduction in the objective function of penalization
terms that are function of higher order moments. The drawback of this
technique is that the ideal weights of the penalization terms are often
unknown. Furthermore, in some cases, an excessive number of higher order
moments may be required to adequately capture all the significant aspect
of the uncertainty embedded into a given problem.  Finally, a wrong choice
of the penalties may lead to a problem formulation that does not have any
feasible solution.  Instead of penalization terms, explicit constraints
can be introduced in the robust optimization problem, and the same
considerations apply for the advantages an the drawbacks of the technique.

Another possibility is the minimax criterion, very popular in statistical
decision theory, according to which the worst case due by uncertainty
is assumed as objective for the optimization. This ensures protection
against worst case scenario, but it is often excessively conservative.

The multi-objective approach \cite{poloni2004robust} based on constrained
optimization is also widely adopted. Here different statistical moments
are used as independent trade-off objectives. The obtained Pareto
front allows an a-posteriori choice of the optimal design between
a set of equally ranked candidates. In this case a challenge is
posed by the increase in the dimensionality of the Pareto front when
several statistical moments are used.  The research related to the
multi-objective method has led to several extensions of the classical
Pareto front concept.  In \cite{teich2001pareto}, for example, the Pareto
front exploration in presence of uncertainties is faced introducing the
concept of {\it probabilistic dominance}, which is an extension of the
classical {\it Pareto dominance}. While in \cite{hughes2001evolutionary},
a probabilistic ranking and selection mechanism is proposed that
introduces the {\it probability of wrong decision} directly in the
formula for rank computation.

An interesting approach, similar in some aspects to the one here
described, is found in \cite{coelho2011multi} where a quantile based
approach is coupled with the probability of {\it Pareto nondominance}
(already seen in \cite{hughes2001evolutionary}). Here, contrary to the
cited work, the optimization technique introduced relies on direct
estimation of the quantile function obtained through the Empirical
Cumulative Distribution Function.

\section{The Generalized Inverse Distribution Function Method}

In the methodology presented herein, the operator that is used to
eliminate the dependence on random variables is the quantile function
of the objective function to be optimized, calculated in one or more
points of its domain of definition.

Before going into the details of the exposure, the definitions of
\emph{Cumulative Distribution Function} (CDF) and \emph{Generalized
Inverse Distribution Function} (GIDF) that will be used are reported.

The ``cumulative distribution function'' (CDF), or just ``distribution
function'', describes the probability that a real-valued random variable
$Q$ with a given probability distribution will be found at a value less
than or equal to $q$.  Intuitively, it is the ``area so far'' function
of the probability distribution.  The CDF is one of the most precise,
efficient and compact ways to represent information about uncertainty,
and a new CDF based approach to robust optimization is here described.

If the CDF is continuous and strictly monotonic then it is invertible, and
its inverse, called quantile function or inverse distribution function,
returns the value below which random draws from the given distribution
would fall, $s\times 100$ percent of the time.  That is, it returns the
value of $q$ such that
\begin{equation} \label{percent}
F_Q(q) = \Pr(Q \le q) = s
\end{equation}
Hence $ F^{-1}( s ),\; s \in [0,1] $ is the unique real number $q$
such that $ F_Q(q) = s $.

Unfortunately, the distribution does not, in general, have an inverse. If
the probability distribution is discrete rather than continuous then
there may be gaps between values in the domain of its CDF, while, if the
CDF is only weakly monotonic, there may be ``flat spots'' in its range.
In general, in these cases, one may define, for $ s \in [0,1] $, the
``generalized inverse distribution function'' (GIDF)
\[
q^s  = Q(s) = F_Q^{ - 1} \left( s \right) = \inf\left\{ {q \in \mathbb{R}:F\left( q \right) \geq s} \right\}
\]
that returns the minimum value of $s$ for which the previous probability statement (\ref{percent}) holds.
The infimum is used because CDFs are, in general, weakly monotonic and
right-continuous (see \cite{CDF}).

\begin{figure}[h!]
\begin{center}
\includegraphics[trim=0 50 0 100,clip,scale=0.45]{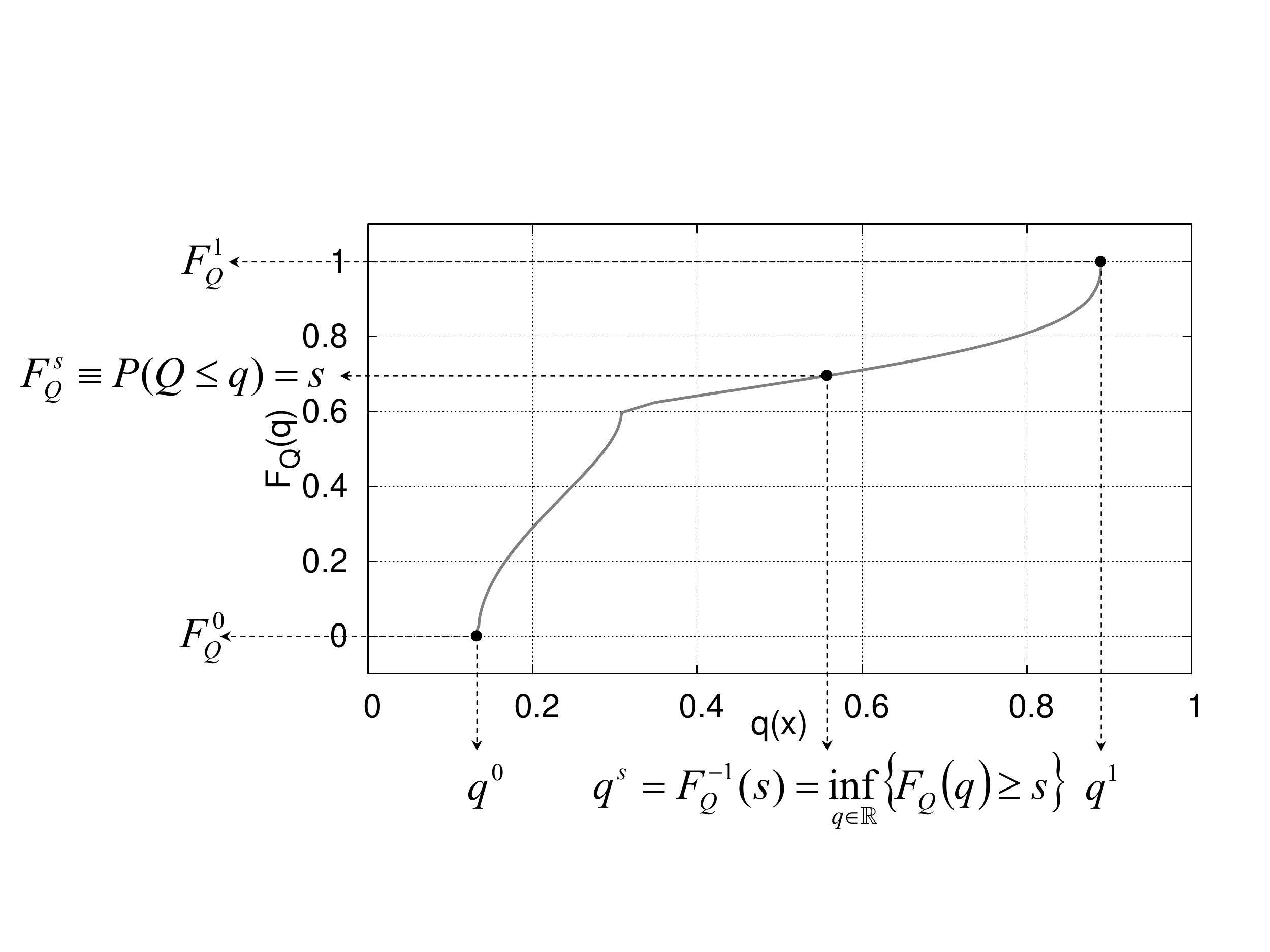}
\end{center}
\caption{CDF and ICDF characteristic points.}
\label{fig02}
\end{figure}

For the sake of clarity, the introduced nomenclature and related
characteristic points $\left( q^s, F_Q^s \right)$ with $q^0\le q^s\le q^1$
and $F_Q^0=0\le F_Q^s\le F_Q^1=1$ are illustrated in Figure \ref{fig02}.

Now that the CDF and the GIDF have been introduced, it becomes immediate
to define, within the framework of multi-objective optimization, a
robust optimization problem in terms of an arbitrary number of quantiles
to optimize:
\begin{equation} \label{quantileopt}
P_{Q(s_i)} : \quad 
\mathop {\min}\limits_{z\in Z} q^{s_i } \left( z \right) = \mathop {\min}\limits_{z\in Z} \inf \left\{ {q\left( z \right) \in \mathbb{R}:F_Q \left( {q\left( z \right)} \right) \geq s_i } \right\}\quad i = 1,\ldots,n
\end{equation}
being $n$ the number of objectives chosen. The approach, then, can be
further extended by introducing objectives that are arbitrary functions
of quantiles.

Of course, the problem now is focused on how to satisfactorily calculate
the quantiles required by the method.  In this work the \emph{Empirical
Cumulative Distribution Function} (ECDF) is used for this purpose. The
definition of ECDF, taken from \cite{ECDF}, is reported here for the
sake of completeness.

Let $X_1,\ldots,X_n$ be random variables with realizations $x_i \in
\mathbb{R}$, the empirical distribution function is an indicator function
that estimates the true underlying CDF of the points in the sample. It
can be defined by using the order statistics $X_{(i)}$ of $X_i$ as:
$$
\widehat{F}_n(x,\omega)=
\begin{cases}
0 & \text{if $x<x_{(1)}$;}\\
\frac{1}{n} & \text{if $x_{(1)}\leq x<x_{(2)}$, $1\leq k<2$;}\\
\frac{2}{n} & \text{if $x_{(2)}\leq x<x_{(3)}$, $2\leq k<3$;}\\
\vdots\\
\frac{i}{n} & \text{if $x_{(i)}\leq x<x_{(i+1)}$, $i\leq k<i+1$;}\\
\vdots\\
1 & \text{if $x\geq x_{(n)}$;}
\end{cases}
$$
where $x_{(i)}$ is the realization of the random variable $X_{(i)}$
with outcome (elementary event) $\omega \in \Omega$.

From now on, therefore, when the optimization algorithm will require
the calculation of the $F_Q(s)$, it will used instead its estimator
$\widehat{F}_{Q_n}(s)$, where $ n $ indicates the number of samples
used to estimate this ECDF.

Note that each indicator function, and hence the ECDF, is itself a
random variable. This is a very delicate issue to consider. Indeed,
if the EDCF is used to approximate the deterministic operator $Q(s)$,
a direct  residual influence of the random variables that characterize
the system under investigation remains on $P_{Q(s)}$. In other words
$Q(s)$ behaves as a random variable, but with the important difference
that its variance tends to zero when the ECDF approximates the CDF
with increasing precision.  It is possible to demonstrate that the
estimator $\widehat{F}_{Q_n}(s)$ is consistent, as it converges almost
surely to $F_Q(s)$ as $n \rightarrow \infty$, for every value of $s$
\cite{Vaart1998}.  Furthermore, for the Glivenko-Cantelli theorem
\cite{Serfling2008}, the convergence is also uniform over $s$.
This implies that, if the ECDF is calculated with sufficient accuracy,
it can be considered and treated as a deterministic operator. On the
other hand, if the number of samples, or the estimation technique of
the ECDF, do not allow to consider it as such, one can still correlate
the variance of the ECDF with the precision of the obtained estimate.
Of course, if the ECDF is estimated in a very precise way, it is possible
to use for the optimization also an algorithm conceived for deterministic
problems, provided that it has a certain resistance to noise. Conversely,
if the ECDF is obtained from a coarse sample, its practical use is only
possible with optimization algorithms specifically designed.

For the same reason, it is often convenient, especially in applications
where the ECDF is defined with few samples, to use $q^\epsilon$ instead
of $q^0$, with $\epsilon>0$ and small, but such that a not excessive
variance of the estimate of $q^\epsilon$ is ensured.

\section{Illustrative example}\label{illexam}

\begin{figure}[h!]
\centerline{
\includegraphics[scale=0.5]{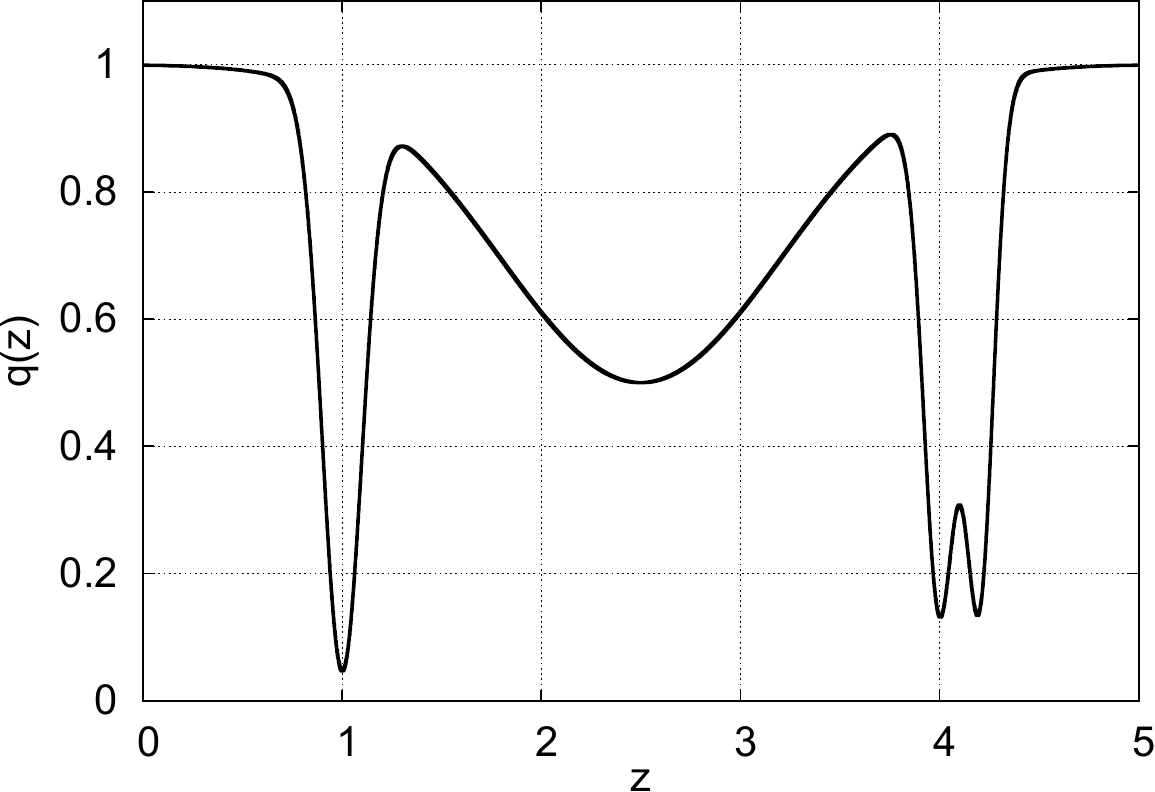}
\includegraphics[scale=0.5]{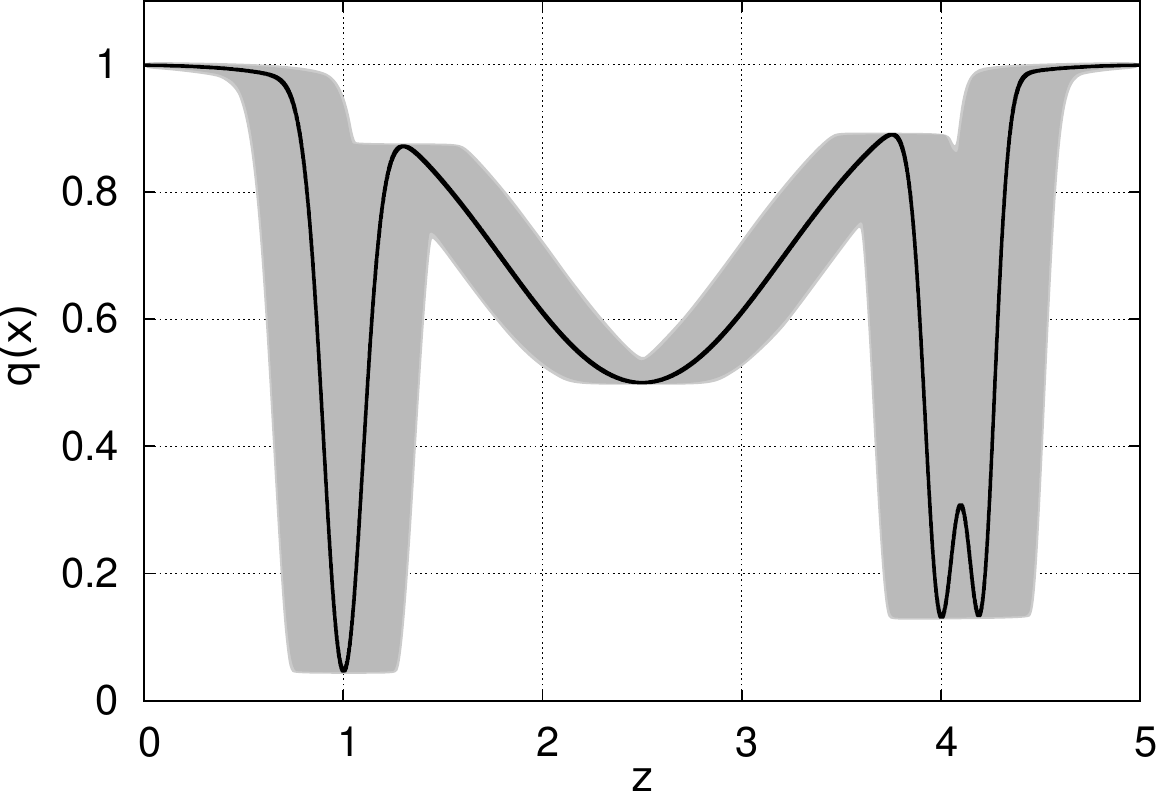}
}
\caption{$q(\mathbf{z},\mathbf{u}=0)$ plot (left) and $q(\mathbf{x})$ projection in the plane $(\mathbf{z},q)$ (right).}
\label{func}
\end{figure}

The features of the quantile curve optimization approach will be
illustrated with the help of a simple example function defined as follows:
\[
q(\mathbf{z},\mathbf{u}) = 1 - \sum\limits_{i = 1}^m {a_i } e^{ - \beta _i \sum\limits_{j = 1}^n {\left( {z_j  + u_j  - c_{i,j} } \right)^2 } } 
\]
with design parameter vector
$
\mathbf{z} = (z_1 , \ldots ,z_n ) \in Z \subseteq \mathbb{R} ^n 
$
and uncertainty vector
$
\mathbf{u} = (u_1 , \ldots ,u_n ) \in U \subseteq \mathbb{R} ^n 
$.
The random variables $\mathbf{u}$ are assumed uniformly distributed
with expected value $0$ and variance $1/12$.  Let's assume for the sake
of compactness:
$
\mathbf{x} = \mathbf{z} + \mathbf{u} = (x_1 , \ldots ,x_n ) = (z_1  + u_1 , \ldots ,z_n  + u_n )
$,
so the example function becomes:
\begin{equation} \label{es}
q(\mathbf{x}) = 1 - \sum\limits_{i = 1}^m {a_i } e^{ - \beta _i \sum\limits_{j = 1}^n {\left( {x_j  - c_{i,j} } \right)^2 } } 
\end{equation}
where the uncertain parameter $\mathbf{u}$ is incorporated in the new
random variable $\mathbf{x}$ Choosing the following set of parameters:
\[
\left\{ {\begin{array}{lll}
   {n = 1} \hfill  \\
   {m = 4} \hfill  \\
   z_1  \in \left[ {0,5} \right],&\multicolumn{2}{c}{u_1  \in \left[ { - 0.25,0.25} \right]} \hfill  \\
   a_1  = 0.9,&\beta _1  = 50.0,&{c}_{{1,1}}  = 1.0 \hfill  \\
   a_2  = 0.5,&\beta _2  = 1.0, &{c}_{{2,1}}  = 2.5 \hfill  \\
   a_3  = 0.8,&\beta _3  = 80.0,&{c}_{{3,1}}  = 4.0 \hfill  \\
   a_4  = 0.8,&\beta _4  = 100.0&{c}_{{4,1}}  = 4.2 \hfill
 \end{array} } \right.
\]
the function reported in Figure \ref{func} is obtained, where
the uncertain parameter $\mathbf{u}$ is incorporated into the
new random variable $\mathbf{x}$ that has expected value equal
to $\mathbf{z}$ and variance equal to $1/12$.  In particular, the
first graph is related to the ``deterministic'' version of function,
e.g.  with $\sigma_\mathbf{u}=0$, while the second one reports the
projection of $q(\mathbf{x})=q(\mathbf{z},\mathbf{u})$ in the plane
$(\mathbf{z},q)$. The black curve corresponds to $ \mathbf{u} = 0 $,
while the contributions to $q$ due to $ \mathbf{u} \neq 0 $ are shown
in gray and indicate the variation of the function $q$ caused by the
random variable $\mathbf{u}$.  Obviously, in the plane $(\mathbf{x},
q)$ the effect of the random variable $\mathbf{u}$ is a simple variation
of the position of the point along the curve $q$, and that is how this
function will be represented in the figures below.

The Empirical CDF is used to estimate
the uncertainty on $q$ induced by $\mathbf{u}$:
\[
\widehat F_Q (q)\approx F_Q (q) = P(Q \leq q)
\]

\begin{figure}[h!]
\centerline{
\includegraphics[scale=0.8]{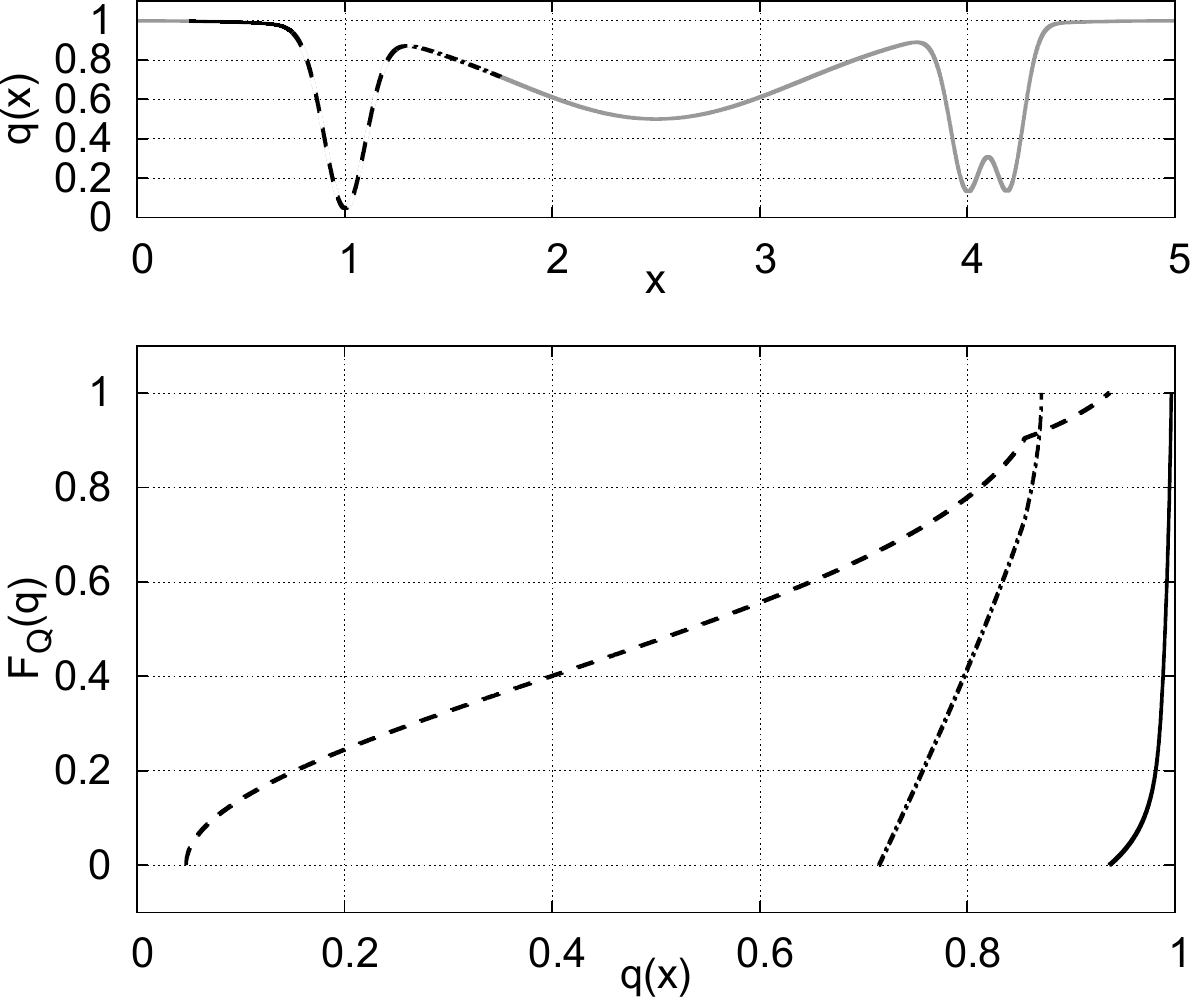}
}
\caption{Some ECDFs related to $q(x)$.}
\label{cdf+f}
\end{figure}

Some ECDFs related to the defined $q(x)$ are reported in Figure
\ref{cdf+f}. An ECDF defined in the $\left(q(x),F_Q(q)\right)$
plane corresponds to a function variation in the $\left(x,q(x)\right)$
plane. The correspondence between ECDF and function variation is evidenced
in the figure by using the same line style.

\begin{figure}[h!]
\centerline{
\includegraphics[scale=0.8]{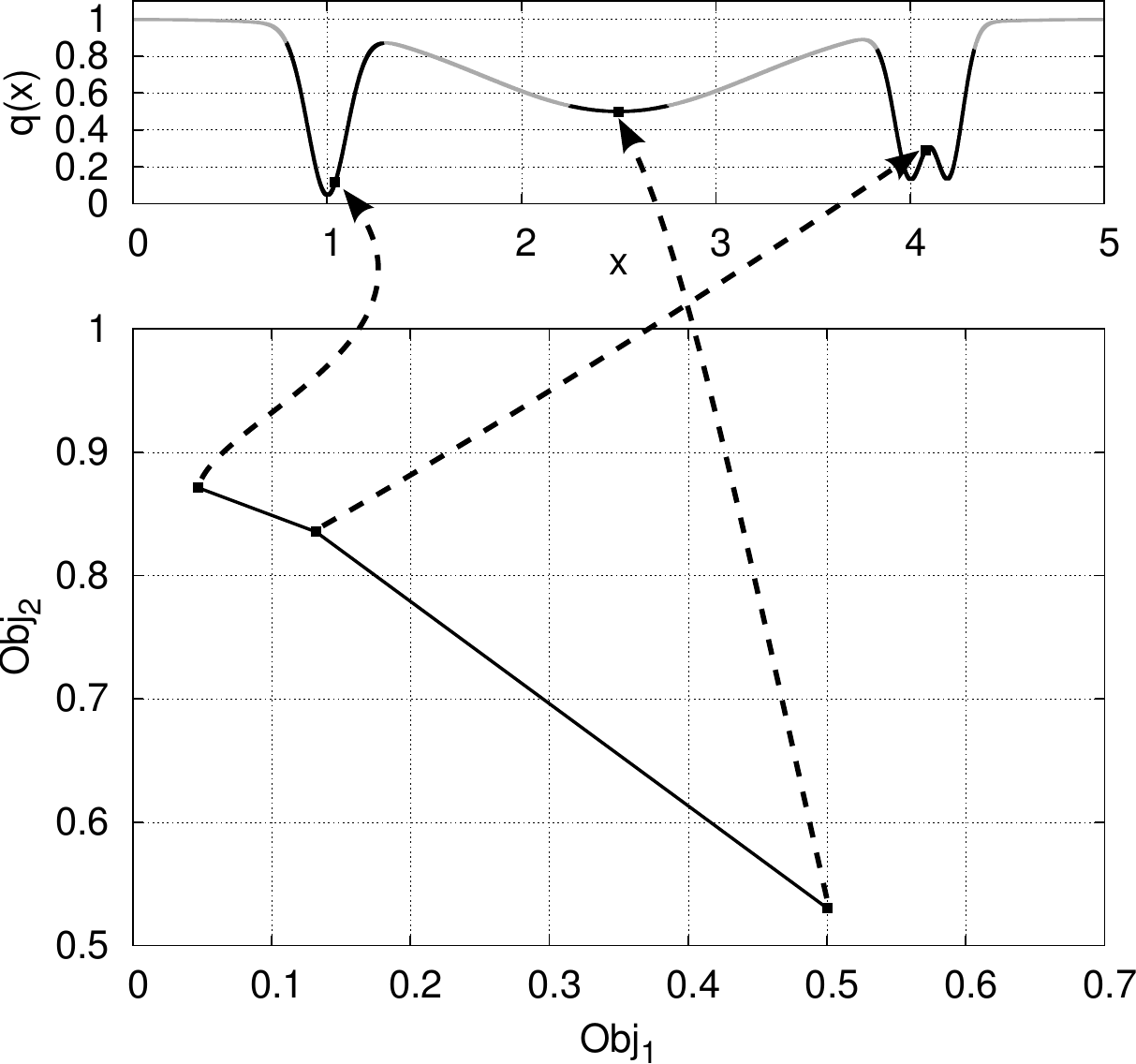}
}
\caption{Pareto front of problem formulation 1.}
\label{pf1}
\end{figure}

Let's now consider the following two-objective problem:
\begin{equation}\label{formulation1}
\mathop {\min}\limits_x \left( {q^\epsilon  ,q^{1-\epsilon} } \right)
\end{equation}
where the use of a small $\epsilon$ value is introduced to account for the
approximation introduced by the ECDF estimator.  Figure \ref{pf1} reports
the results obtained with problem formulation (\ref{formulation1}). The
extreme of the Pareto front related to the best $\mathrm{Obj}_1$,
namely $q^\epsilon$, is representative of the best possible optimum,
without regard to the variance. The front extreme related to the best
$\mathrm{Obj}_2$ gives, instead the most robust solution, e.g. the one
that that has the least variance.  Finally, the solution located in
the middle of the front represents a compromise between best absolute
performance and smallest variance.

\begin{figure}[h!]
\centerline{
\includegraphics[scale=0.8]{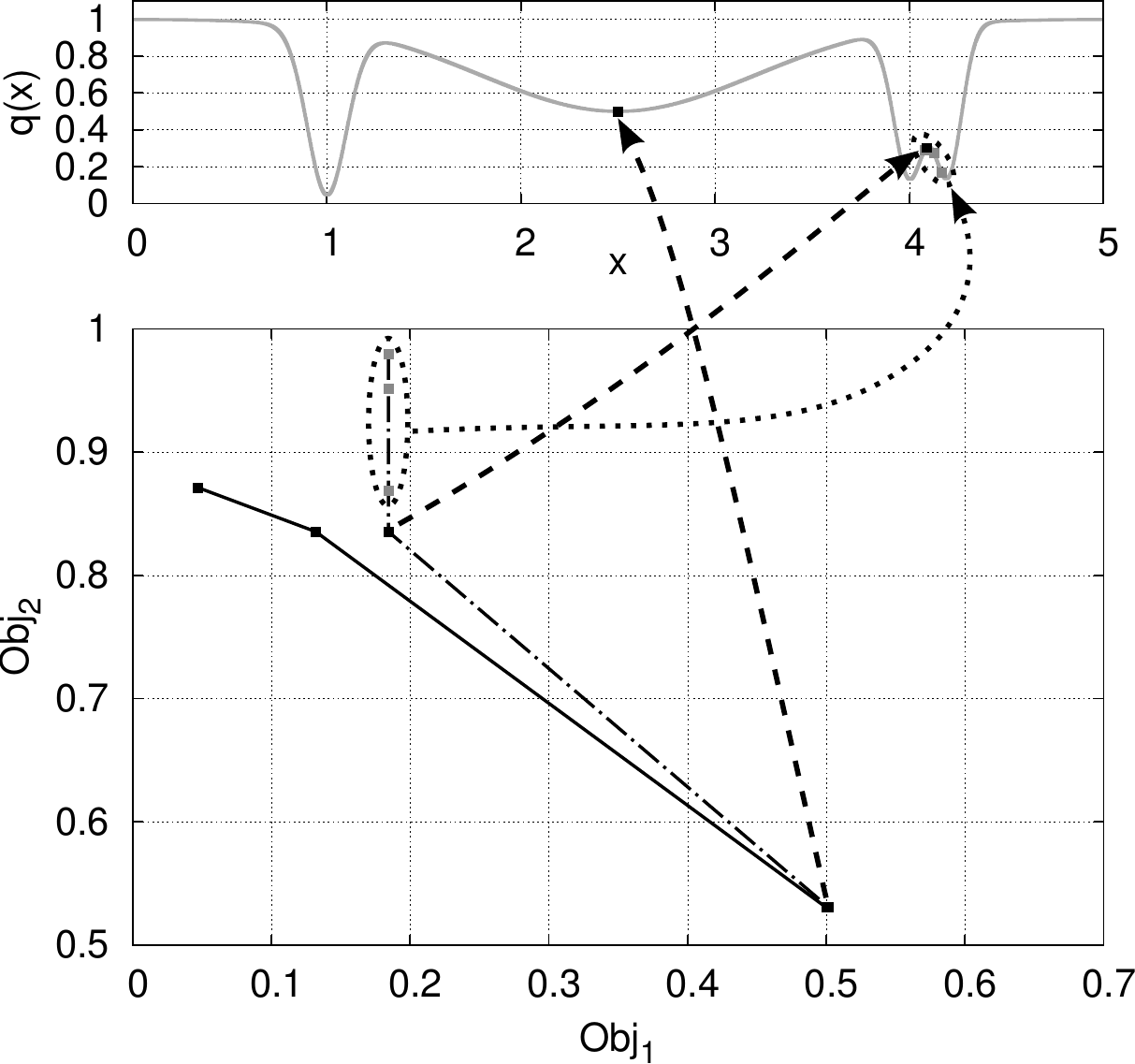}
}
\caption{Pareto front of problem formulation 2.}
\label{pf2}
\end{figure}

Figure \ref{pf2} reports, as a dot-dashed line, the Pareto front obtained
solving the problem
\begin{equation}\label{formulation2}
\mathop {\min}\limits_x \left( {q^{0.25} ,q^{1-\epsilon} } \right)
\end{equation}
Here the solutions with higher variance are no longer present on the
Pareto front, and this is often a desirable behavior in a robust design
problem. Among the solutions which are on the left side of the front,
the one with lowest $\mathrm{Obj}_2$ should be selected. Indeed, the
significant worsening observed in the second objective is not compensated
by an adequate improvement of the first one.

\begin{figure}[h!]
\centerline{
\includegraphics[scale=0.8]{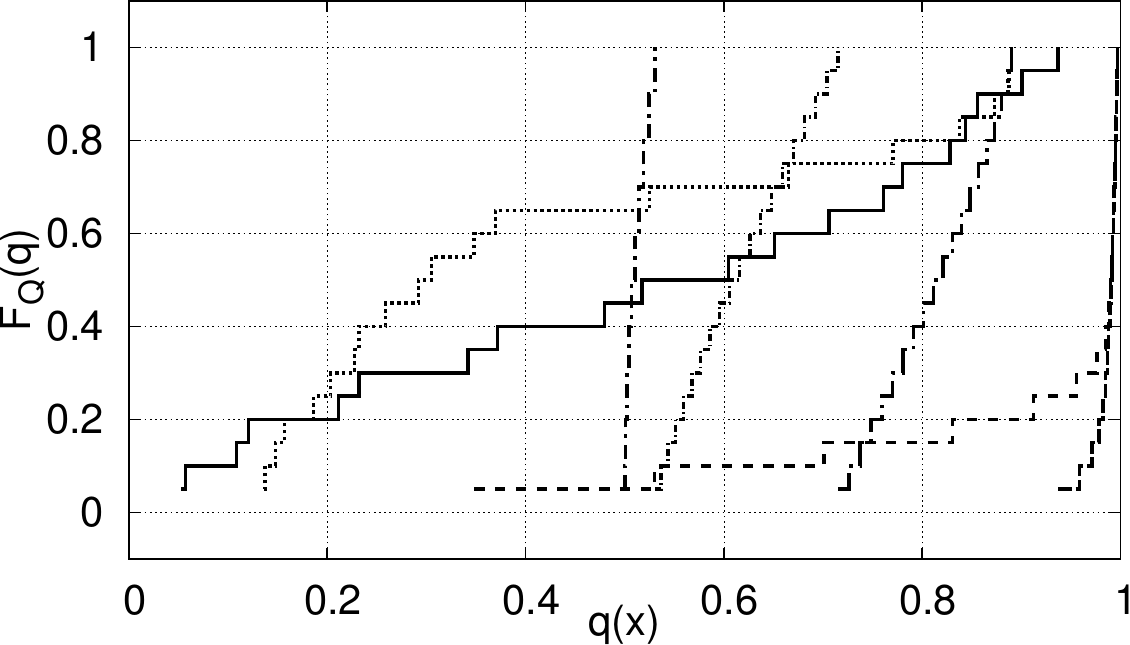}
}
\caption{ECDF estimated using just 20 uniform samples.}
\label{ecdf20}
\end{figure}

\begin{figure}[h!]
\centerline{
\includegraphics[scale=0.8]{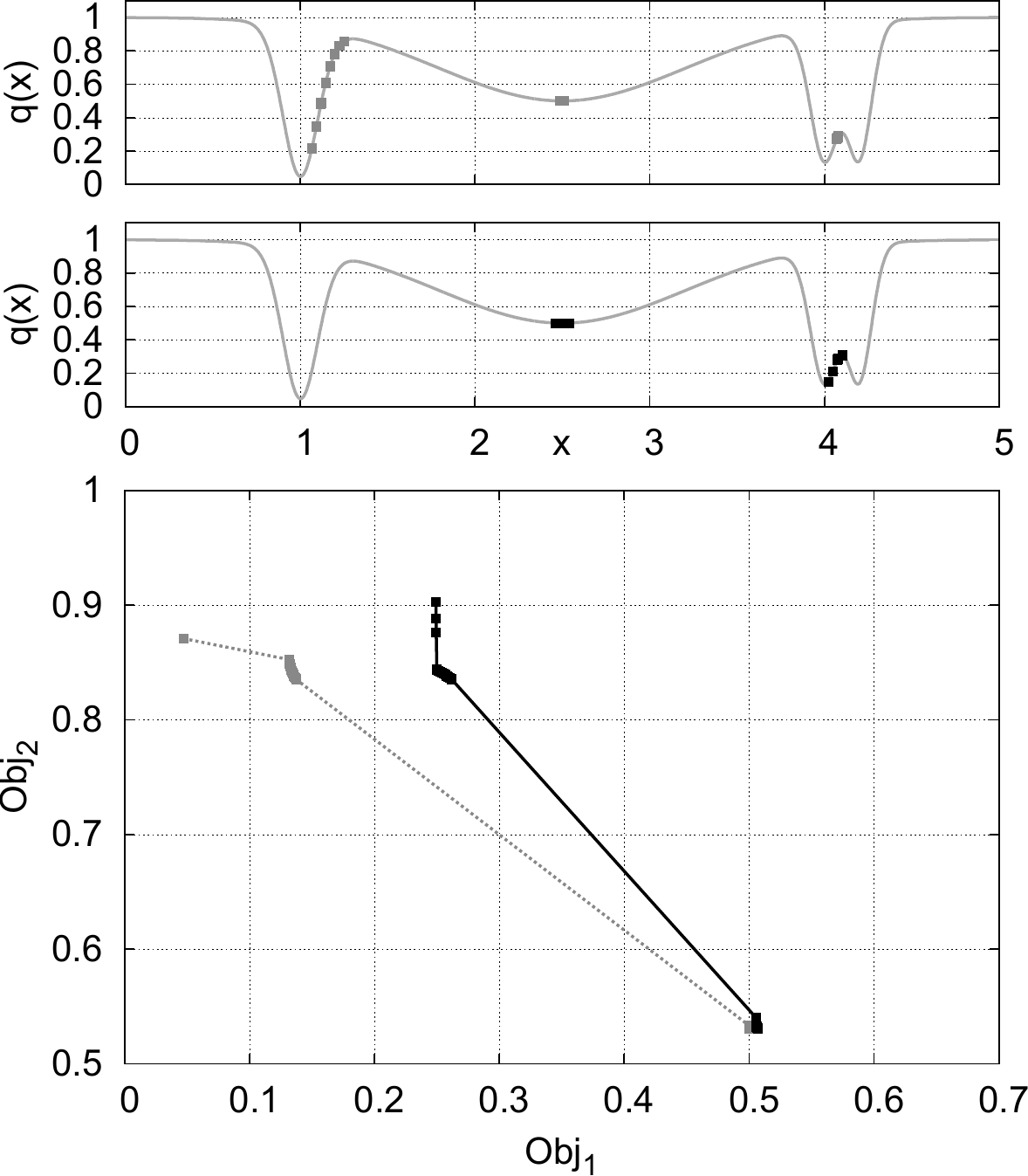}
}
\caption{Comparisons of Pareto fronts and related results of problems $\left( {q^{\epsilon} ,q^1 } \right)$, dotted curve, and $\left( {q^{0.45} ,q^1 } \right)$, solid curve.}
\label{pf3}
\end{figure}

Finally, it is worth to see what should be expected when the estimate of
the CDF is very coarse. Consider, in this regard, a ECDF calculated with
only 20 samples distributed uniformly (see Figure \ref{ecdf20}). Here,
both the problem (\ref{formulation1}) and problem (\ref{formulation2})
do not allow to discern between the solution with maximum variance and
the one giving an intermediate compromise. If, instead, the following
problem is solved
\begin{equation}\label{formulation3}
\mathop {\min}\limits_X \left( {q^{0.45} ,q^{1-\epsilon} } \right)
\end{equation}
it is possible to exclude the solution with maximum variance from the
Pareto front, as can be observed by analyzing Figure \ref{pf3}. Thus,
this approach is definitely worth to be tested when a sufficiently
precise estimate of the CDF is not available.

\section{A benchmark function with variable number of design parameters}

The function reported in Table \ref{tab:1}, taken from \cite{Vasile2011},
is used as a benchmark to test the GIDF based approach to robust
optimization.  With respect to the function reported in the reference,
the following changes have been introduced: the ranges of design and
uncertain parameters have been changed as reported in table, and a
multiplicative factor equal to $1/n$ has been introduced to make easier
the result comparison when the dimension of the parameter space changes.
The random variables $\mathbf{u}$ have a uniform distribution function.
Table \ref{tab:2} reports the solutions to the optimization problems
$$\min\limits_{\mathbf{d}\in D, \mathbf{u}\in U} f(\mathbf{d},\mathbf{u})$$
$$\min\limits_{\mathbf{d}\in D} \max\limits_{\mathbf{u}\in U} f(\mathbf{d},\mathbf{u})$$
over the cartesian product of $D$ and $U$. The first problem represents
the best possible solution obtainable if the $\mathbf{u}$ are considered
as design parameters varying in $U$. The second one, instead, minimizes
the maximum possible loss or, alternatively, maximizes the minimum
gain, according to the framework of decision theory \cite{Neumann1953}.
These solutions have been obtained analytically and verified by exhaustive
search for $n=1$. It is worth to note that these particular optimal
solutions are the same whatever is the dimension of the search space.

The optimization algorithm used here is a simple multi-objective genetic
algorithm not specially conceived for optimization under uncertainty.
The algorithm is based on the Pareto dominance concept \cite{Deb:2001}
and on local random walk selection \cite{Quaglia:2000b,Vicini:97b}.
Crossover operator is the classical \emph{one-point} crossover which
operates at bit level, while mutation operator works at the level of
the design vector parameters (which are real numbers). A parameter,
called \emph{mutation rate} controls the operator activation probability
for each variable vector element, while a further parameter, called
\emph{strength}, is the maximum relative value the uniform word mutation.
The word mutation value is given by $\mathrm{strength}\cdot(r-0.5)(u-l)$
with $r \in [0,1]$ uniform random number, $u$ upper variable bound and $l$
lower variable bound. An elitist strategy was adopted in the optimization
runs. It consists in replacing 20\% of the population calculated at
each generation with elements taken at random from the current Pareto
front. Obviously, the elements of the population are used to update the
current Pareto front before the replacement, in order to avoid losing
non-dominated population elements.

The multi-objective runs were performed using 100\% crossover activation
probability and word mutation with \emph{mutation rate} equal to 50\% and
\emph{strength} equal to 0.06.  The initial population was obtained using
the quasi-random low-discrepancy Sobol sequence \cite{Bratley1988}.  The
ECDF used to estimate the CDF was obtained with 2500 Montecarlo samples
in all runs. The population size was set to 4000 elements for all runs,
while the number of generations was set to $10$ for $n=1$, $200$ for $n=2$
and $1000$ for $n=6$.  The problem solved was $\min\limits_{\mathbf{z}\in
Z} \left( {q^\epsilon  ,q^{1-\epsilon} } \right)$.

\begin{table}[h!]
\caption{Benchmark functions table.}
\label{tab:1}
{\scriptsize
\begin{center}
\begin{tabular}{cccc}
\hline\noalign{\smallskip}
ID &	Function &	Ranges &	Dimension \\
\noalign{\smallskip}\hline\noalign{\smallskip}
MV4  & $
f = \frac{1}{n}\sum\limits_{i = 1}^n {\left( {2\pi  - u_i } \right)\cos\left( {u_i  - d_i } \right)} 
$
&
$
\mathbf{u} \in \left[ {0,3} \right]^n ,\mathbf{d} \in \left[ {0,2\pi } \right]^n 
$
&
	 	 	1, 2 and 6\\
\noalign{\smallskip}\hline\noalign{\smallskip}
\end{tabular}
\end{center}
}
\end{table}

\begin{table}[h!]
\caption{Benchmark functions table results.}
\label{tab:2}
{\scriptsize
\begin{center}
\begin{tabular}{ccccccc}
\hline\noalign{\smallskip}
ID &
\multicolumn{3}{c}{
$\min\limits_{\mathbf{d}\in D, \mathbf{u}\in U} f(\mathbf{d},\mathbf{u})$ 
}
&
\multicolumn{3}{c}{
$\min\limits_{\mathbf{d}\in D} \max\limits_{\mathbf{u}\in U} f(\mathbf{d},\mathbf{u})$
} \\
\noalign{\smallskip}\hline\noalign{\smallskip}
 & $\qquad\mathbf{d}\qquad$ & $\qquad\mathbf{u}\qquad$ & $f$ &
   $\qquad\mathbf{d}\qquad$ & $\qquad\mathbf{u}\qquad$ & $f$ \\
\noalign{\smallskip}\hline\noalign{\smallskip} 
MV4& $[3.1416]^n$ & $[0]^n$ & $ -6.283185\ldots$ &
     $[4.6638]^n$ & $[0]^n$ & $ -0.305173\ldots$ \\
\noalign{\smallskip}\hline\noalign{\smallskip}
\end{tabular}
\end{center}
}
\end{table}

\begin{figure}[h!]
\centerline{
\includegraphics[scale=0.8]{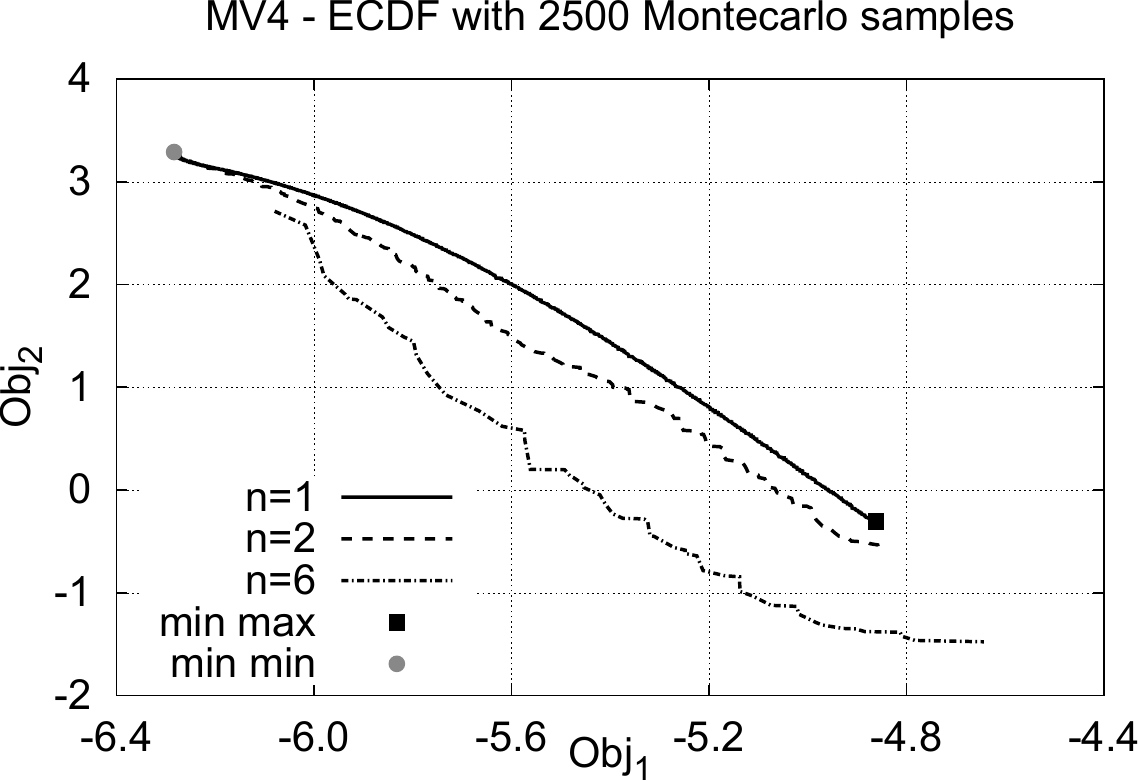}
}
\caption{Pareto fronts and deterministic $\min$ and $\min \max$ solutions for the MV4 test case.}
\label{mv4pf}
\end{figure}

\begin{figure}[h!]
\centerline{
\includegraphics[scale=0.8]{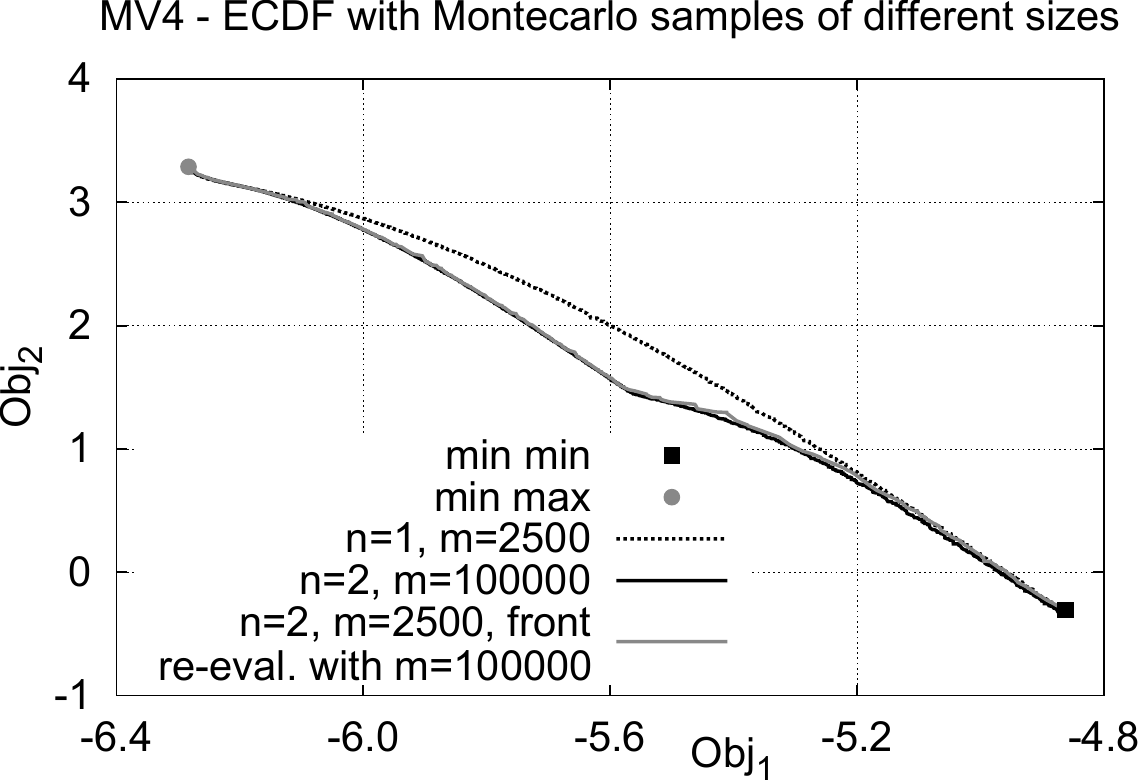}
}
\caption{Pareto fronts for the MV4 test case obtained with different sizes for Montecarlo sampling.}
\label{mv4pfreeval}
\end{figure}

Figure \ref{mv4pf} reports the Pareto fronts and the deterministic
$\min$ and $\min \max$ solutions obtained for the MV4 test case at
different values of the design space size $n$.  It can be easily observed
that, in the case $n = 1$, the extremes of the front are practically
coincident with the deterministic solutions, while, in the case $n = 2$,
the solution of the Pareto front which minimizes the second objective
$(q^{1-\epsilon})$ underestimates the $\min\max$ solution. The trend
is even more evident in the case $n = 6$, where, indeed, also the
extreme of the front that minimizes the first goal $(q^{\epsilon})$
overestimates the value obtained from the $\min$ problem.  This can be
explained by the fact that the two deterministic solutions are located in
correspondence with the extremes of variation of the random variables of
the problem. Therefore, as the number of random variables increases, in
accordance with the central limit theorem \cite{Sobol1994}, it becomes
less likely that all random variables are located in correspondence
of one of their limits of variation.  Anyway, as illustrated in Figure
\ref{mv4pfreeval}, when the Pareto front obtained with the sample size
$m$ equal to 2500 is re-evaluated with a larger Montecarlo sample, the
obtained curve is a quite acceptable approximation of the Pareto front
obtained with $m=100000$.

\begin{figure}[h!]
\centerline{
\includegraphics[scale=0.8]{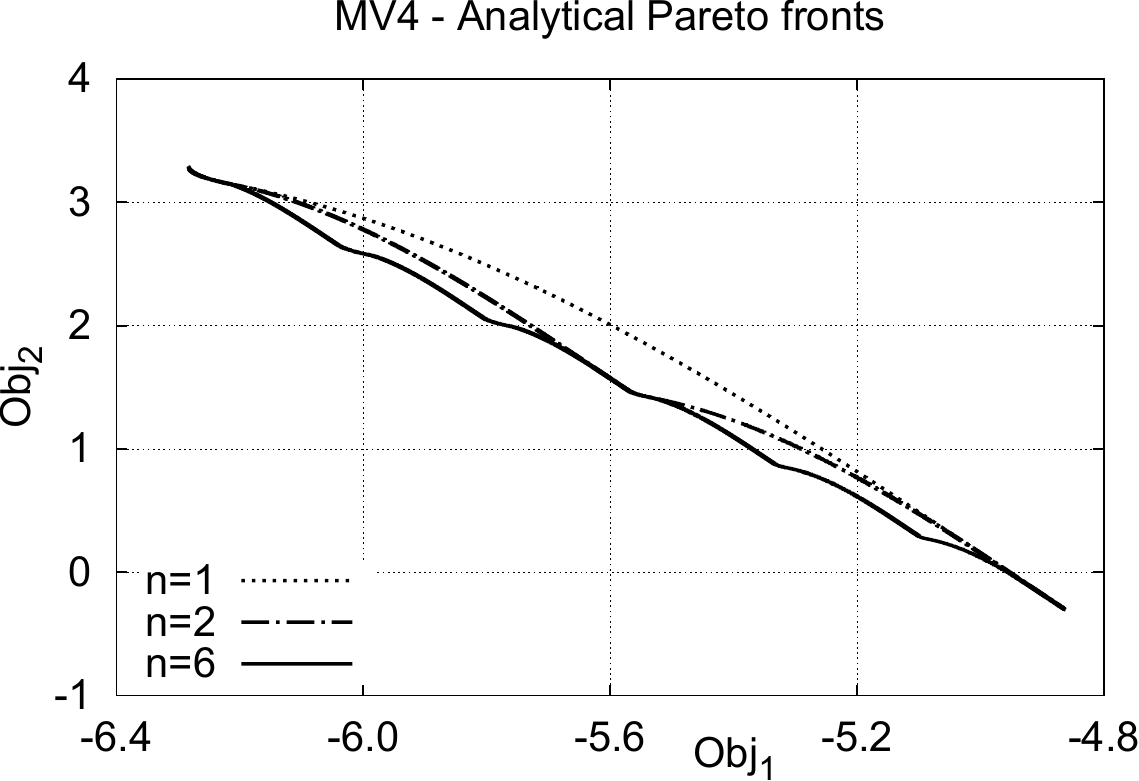}
}
\caption{Pareto fronts for the MV4 test case obtained with the exact method.}
\label{mv4pfexact}
\end{figure}

It should be noted that, for the problem at hand, it is possible to determine a close approximation of the exact Pareto front for any value of $ n $, once a good approximation is known (or the exact Pareto front) the for $ n = 1 $ case.
Note, in this regard, that if $f (\mathbf{x})$, with $n \neq 1$
and $\mathbf{x} = \left(x_1,\ldots,x_n\right)$, is non-dominated, then so are all its components $f_1 (x_1), \ldots, f_1 (x_n)$.
So the elements belonging to the Pareto front with $n \neq 1$ may be constructed by adding elements belonging to the front with $n = 1$. However, not all these elements belong to the Pareto front (with $n\neq 1$), for which it will be necessary to extract the non-dominated elements  among all those obtained in this way.
Figure \ref{mv4pfexact} reports the Pareto fronts at $n=2$ and $n=6$ obtained with this technique.

\begin{figure}[h!]
\centerline{
\includegraphics[scale=0.8]{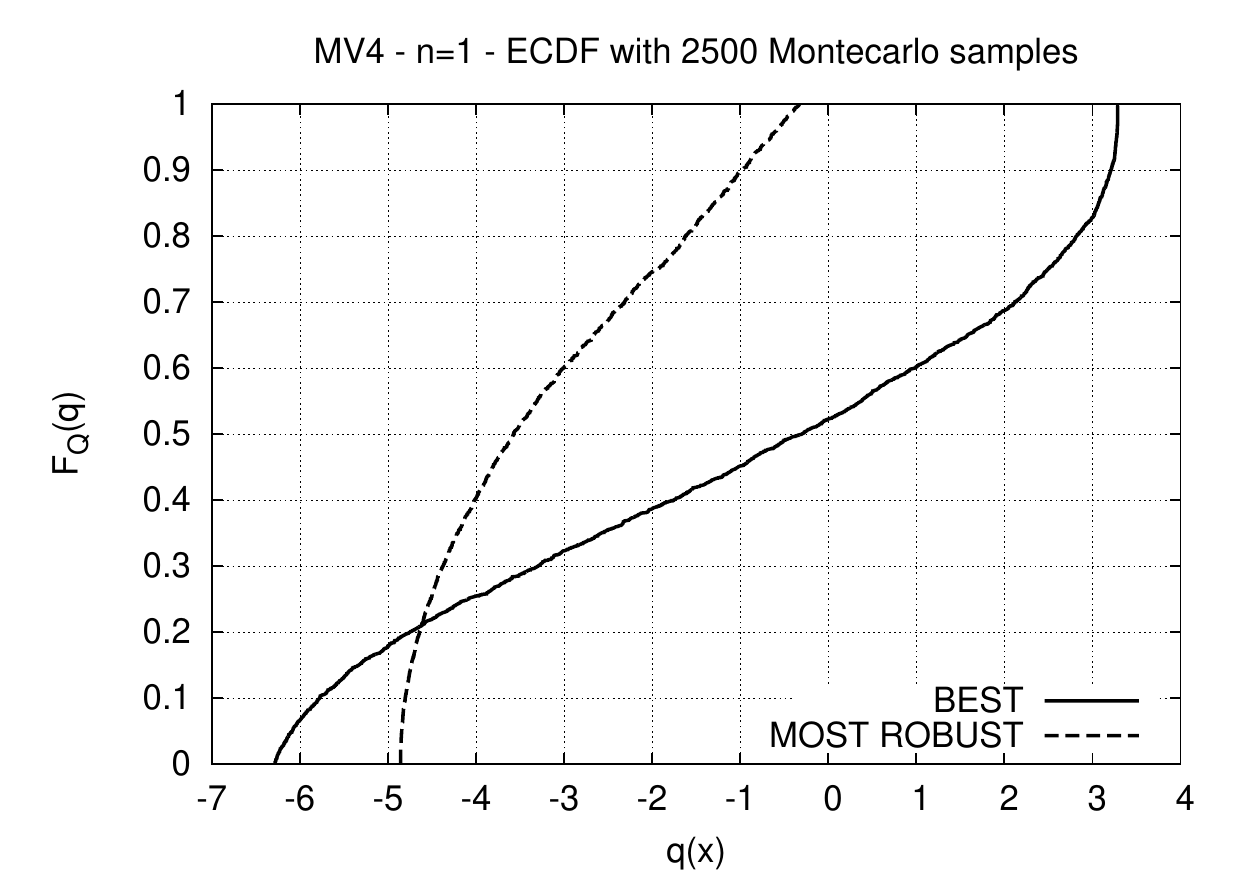}
}
\caption{Optimal ECDF curves for the MV4 with $n=1$.}
\label{mv4cdf}
\end{figure}

\begin{figure}[h!]
\centerline{
\includegraphics[scale=0.8]{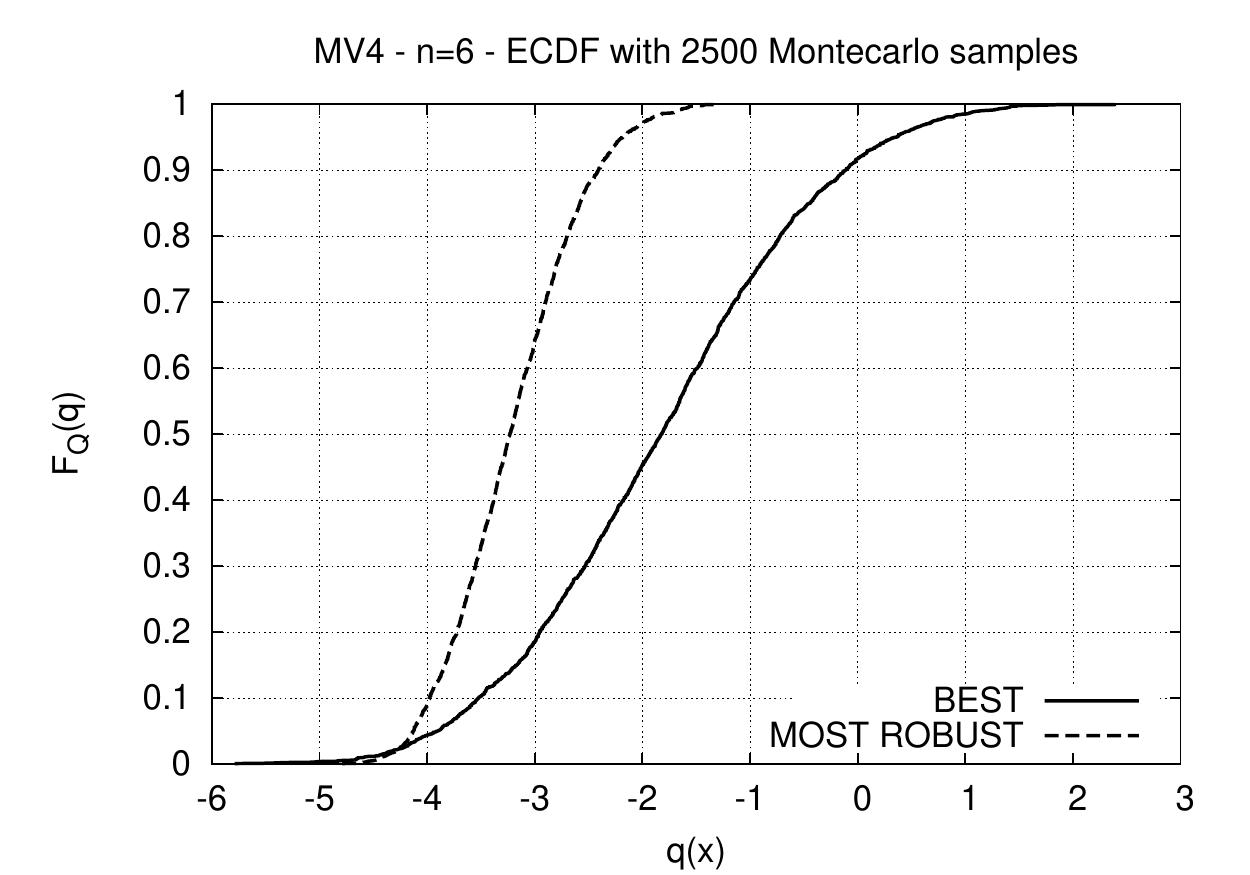}
}
\caption{Optimal ECDF curves for the MV4 with $n=6$.}
\label{mv4cdfn6}
\end{figure}

Figures \ref{mv4cdf} and \ref{mv4cdfn6} show the ECDF corresponding to
the extremes of the Pareto front, respectively for the cases $n = 1$
and $n = 6$. It is noted, again in accordance with the central limit
theorem, that, in the case $n = 6$, the ECDF curves are very close to
those related to a Gaussian distribution.

\section{Evaluating and improving the quantile estimation}

The example in the previous section shows very clearly that the results
of the proposed method may depend in an essential way on the quality
of estimation of quantiles that is obtained through the ECDF. This
leads in a natural way to deal with two issues: how to evaluate the
quality of the estimation of the quantiles used in the multi-objective
optimization problem, and how to possibly get a better quantile estimate
with a given computational effort.  Related to these two points, however,
there are other problems too: it is possible to conceive an algorithm
that can account for the error in the estimation of the quantiles? is it
possible to find estimators of ECDF that give the required accuracy in
the computation of the quantiles of interest for the optimization problem?

The approach here proposed for assessing the quality of the estimate
of the quantile is based on the bootstrap method introduced by Efron in
1977 \cite{Efron1979,Diaconis1983}.

This method represents a major step forward in the statistical
practice because it allows to accurately assess the variability
of any statistical estimator without making any assumption about
the type of distribution function involved.  Suppose that a
statistic $T\left(x_1,x_2,\ldots,x_n\right)$ is given, evaluated
on a set of data $x_1,x_2,\ldots,x_n$ belonging to an assigned
space $X$.  The bootstrap essentially consists in repeatedly
recalculate the statistic $T$ employing a tuple of new samples
$x^*_1,x^*_2,\ldots,x^*_n$ obtained by selecting them from the collection
$\left\{x_1,x_2,\ldots,x_n\right\}$ with replacement. The repeated
calculation of $T\left(x^*_1,x^*_2,\ldots,x^*_n\right)$ gives a set of
values that gives a good indication of the distribution of $T$.

Therefore, to calculate the accuracy of a generic quantile $q^s$, obtained
by the estimator $\widehat{F}_{Q_n}(s)$, the bootstrap procedure can be
applied to the samples that define the estimator. This allows to calculate
the corresponding distribution of $q^s$ for a fixed value of $s$.

\begin{figure}[h!]
\centerline{
\includegraphics[scale=0.8]{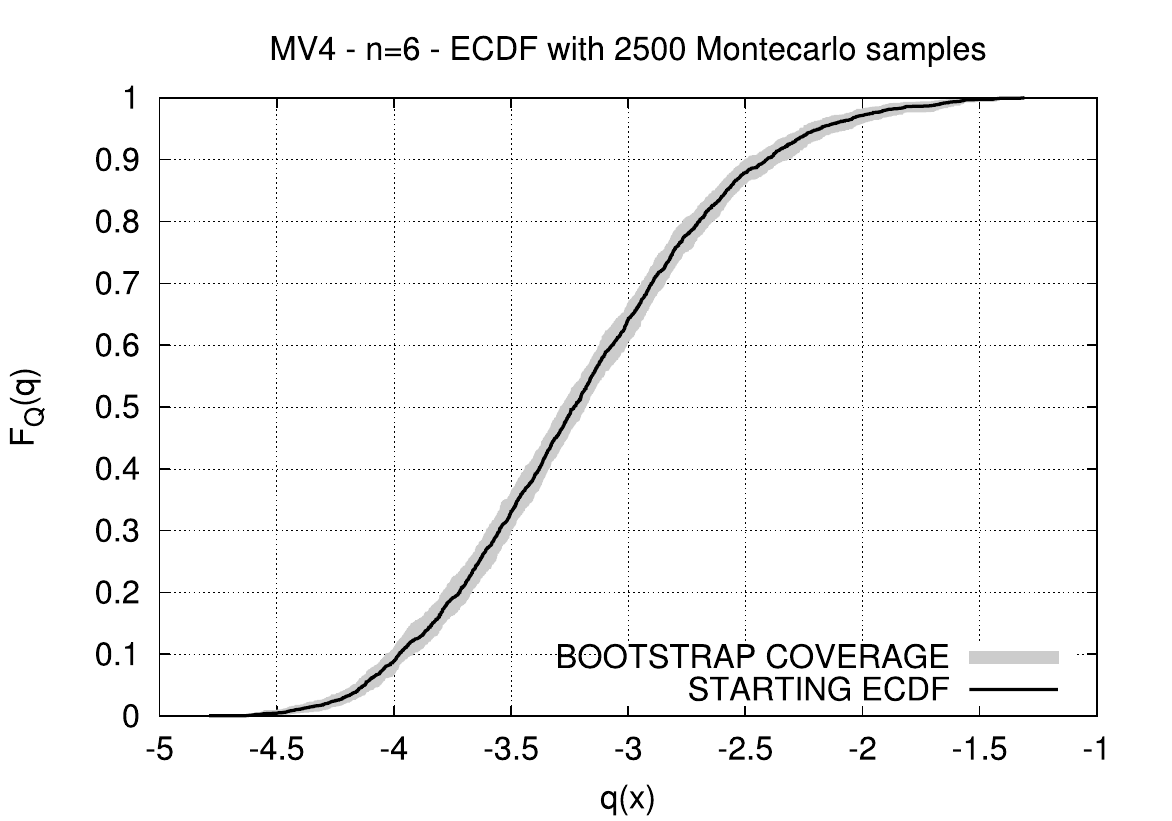}
}
\caption{ECDF corresponding to the most robust solution and related bootstrap coverage.}
\label{bootstrap-coverage}
\end{figure}

\begin{figure}[h!]
\centerline{
\includegraphics[scale=0.8]{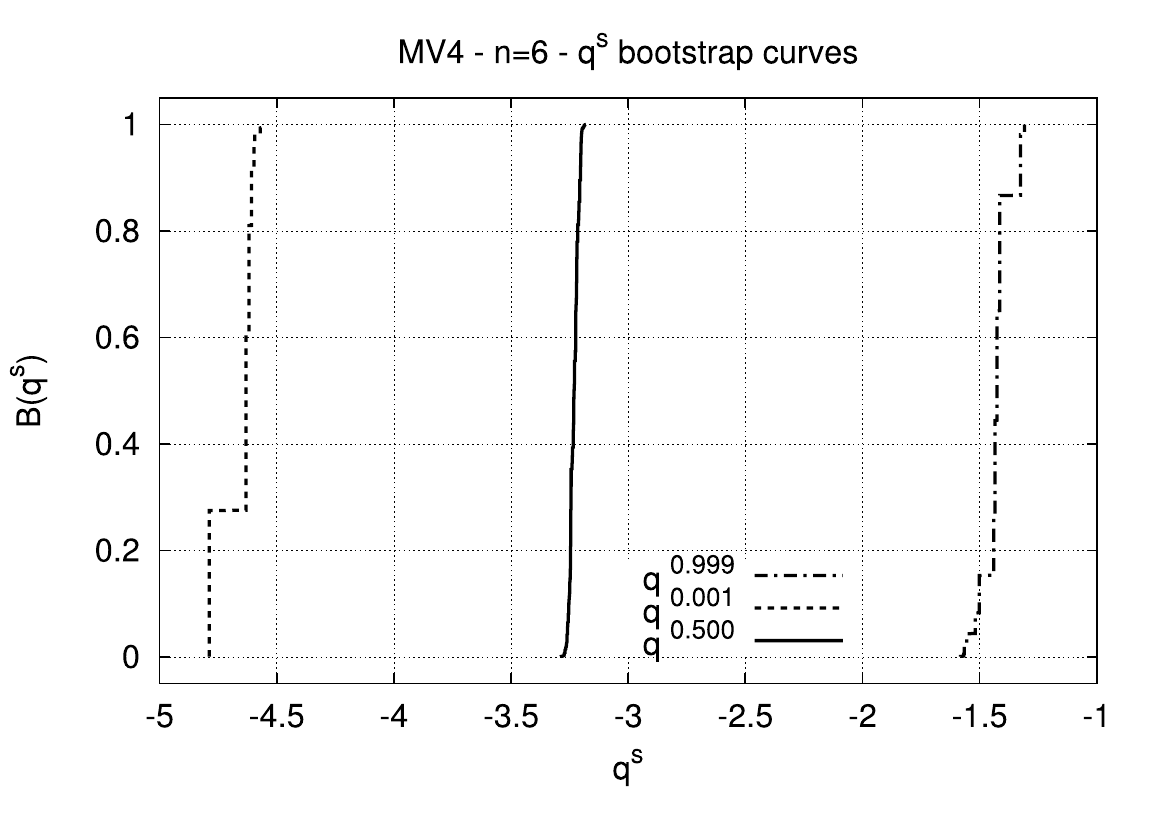}
}
\caption{Cumulative distributon curves accounting for the uncertainty of quantiles evaluation obtained with the bootstrap method.}
\label{bootstrap}
\end{figure}

Figure \ref{bootstrap-coverage} reports the ECDF related to the
solution labeled as ``MOST ROBUST'' in Figure \ref{mv4cdfn6}. The
bootstrap was applied to this ECDF repeating the sampling process
2000 times. The area in gray color represents the superposition
of all the curves obtained in this way. From the bootstrap data
it is then possible to evaluate the accuracy of a given quantile
estimate. Figure \ref{bootstrap} reports the ECDF of three different
quantiles (namely $q^{0.001},q^{0.500},q^{0.999}$) obtained from the
``MOST ROBUST'' solution to MV4.  According to \cite{Diaconis1983},
an accuracy measure for $q^s$ can be obtained considering the central
68\% of bootstrap samples.  These values lay between an interval
$[q^s_\ell,q^s_u]$ centered on the observed value $q^s$. Half the length
of this interval gives a measure of the accuracy for the observed value
that corresponds to the traditional concept of ``standard error''. Here
this value vill be indicated with $\widehat{SE}$ to distinguish
it from the true standard error $SE$.  \footnote{See Wikipedia article
``http://en.wikipedia.org/wiki/Standard\_error\#Standard\_error\_of\_mean\_versus\_standard\_deviation''.
}

Table \ref{accuracy} reports the computed accuracy values for the
considered quantiles for the above cited ``MOST ROBUST'' solution obtained
from an ECDF with 2500 Montecarlo samples. The fourh column reports,
finally, the maximum estimated error $\widehat{ME}$.

\begin{table}
\caption{Quantile estimates and related accuracy for MV4 ``MOST ROBUST'' solution with $n=6$.}
\label{accuracy}
\begin{center}

\begin{tabular}{cccc}
\hline\noalign{\smallskip}
$s$ & $q^s$ & $\widehat{SE}$ & $\widehat{ME}$\\
\noalign{\smallskip}\hline\noalign{\smallskip}
$0.001000$ & $-4.630433$ & $\pm 0.090423$ & $\pm 0.117169$\\
$0.500000$ & $-3.230388$ & $\pm 0.018834$ & $\pm 0.054983$\\
$0.999000$ & $-1.425868$ & $\pm 0.013192$ & $\pm 0.136330$\\
\noalign{\smallskip}\hline\noalign{\smallskip}
\end{tabular}
\end{center}
\end{table}

\begin{figure}[h!]
\centerline{
\includegraphics[scale=0.8]{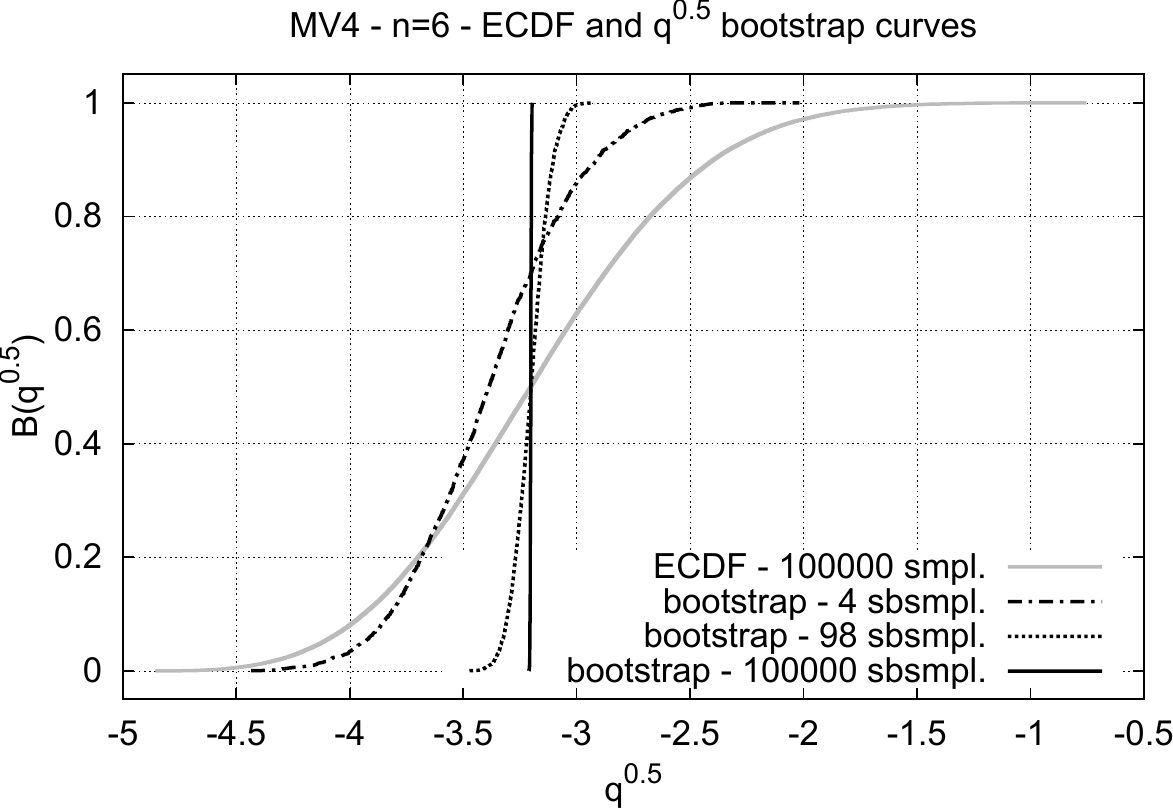}
}
\caption{Bootstrap curves for $q^{0.5}$ and high resolution ECDF for the MV4 ``MOST ROBUST'' solution with $n=6$.}
\label{ECDF+bootstrap-ed}
\end{figure}

\begin{figure}[h!]
\centerline{
\includegraphics[scale=0.8]{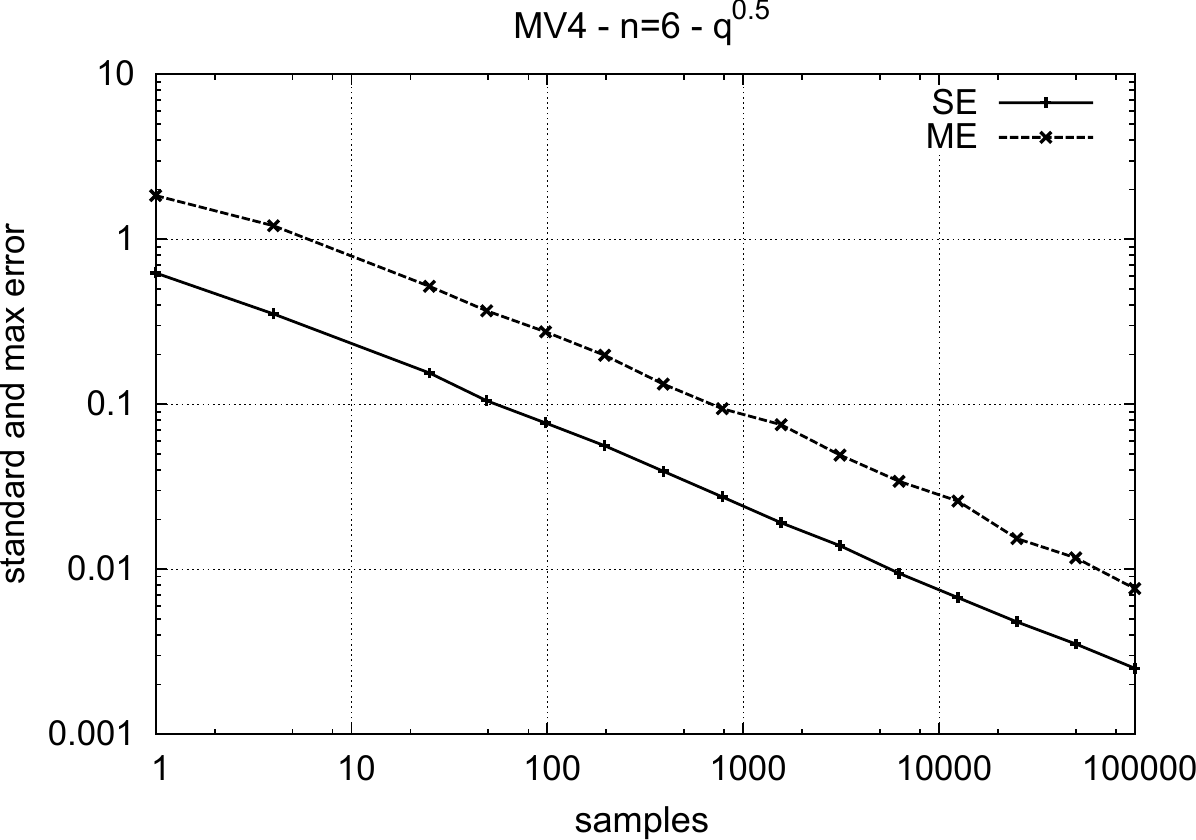}
}
\caption{Standard and maximum error for $q^{0.5}$ as a function of the number of samples for the MV4 ``MOST ROBUST'' solution with $n=6$.}
\label{error-vs-samples-ed}
\end{figure}

The bootstrap is useful for characterizing the standard error of a given
sample used to construct the ECDF, but it is even more important to try
to describe how the error varies depending on the number of elements of
the sample.  To obtain this estimate, a variant of the bootstrap process
is proposed here in which the starting point is a ECDF calculated with
many Montecarlo samples, so as to get as close as possible to the true
CDF. From these data, the bootstrap is applied by drawing a number of
samples smaller than the size of the original sampling. \footnote{for
practical reasons, in the implementation used here, the samples are
replicated so that the resulting ECDF has the same number of points as
the original one.} Note that if the number of samples of the bootstrap
is equal to one, we obtain the starting ECDF.  As can be seen from
Figure \ref{ECDF+bootstrap-ed}, the bootstrap curves\footnote{Maybe
we should call them differently!} tend to widen when the number of
used samples decreases.  From these curves it is possible to derive the
estimates of standard and maximum errors with the same procedure used to
construct Table \ref{accuracy}. Figure \ref{error-vs-samples-ed} shows,
on a logarithmic scale, $\widehat{SE}$ and $\widehat{ME}$ as a function
of the number of samples used to bootstrap $q^{0.5}$.  The important
thing to note is that, already with a number of elements that ranges
from 10 to 100, values of $\widehat{SE}$ and $\widehat{ME}$ are
got that an ad hoc designed optimization algorithm could handle.

\begin{figure}[h!]
\centerline{
\includegraphics[scale=0.8]{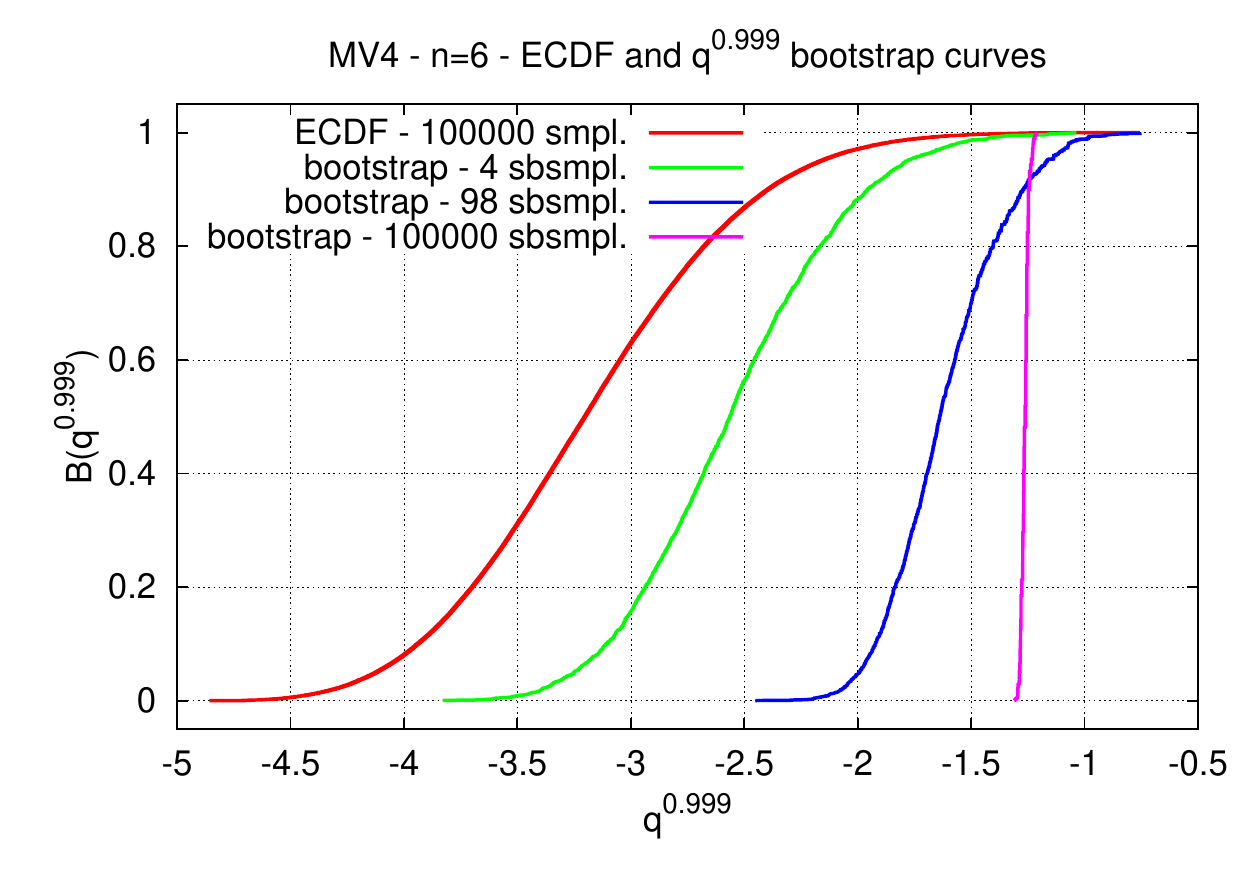}
}
\caption{Bootstrap curves for $q^{0.999}$ and high resolution ECDF for the MV4 ``MOST ROBUST'' solution with $n=6$.}
\label{ECDF+bootstrap-ed-0999}
\end{figure}

\begin{figure}[h!]
\centerline{
\includegraphics[scale=0.8]{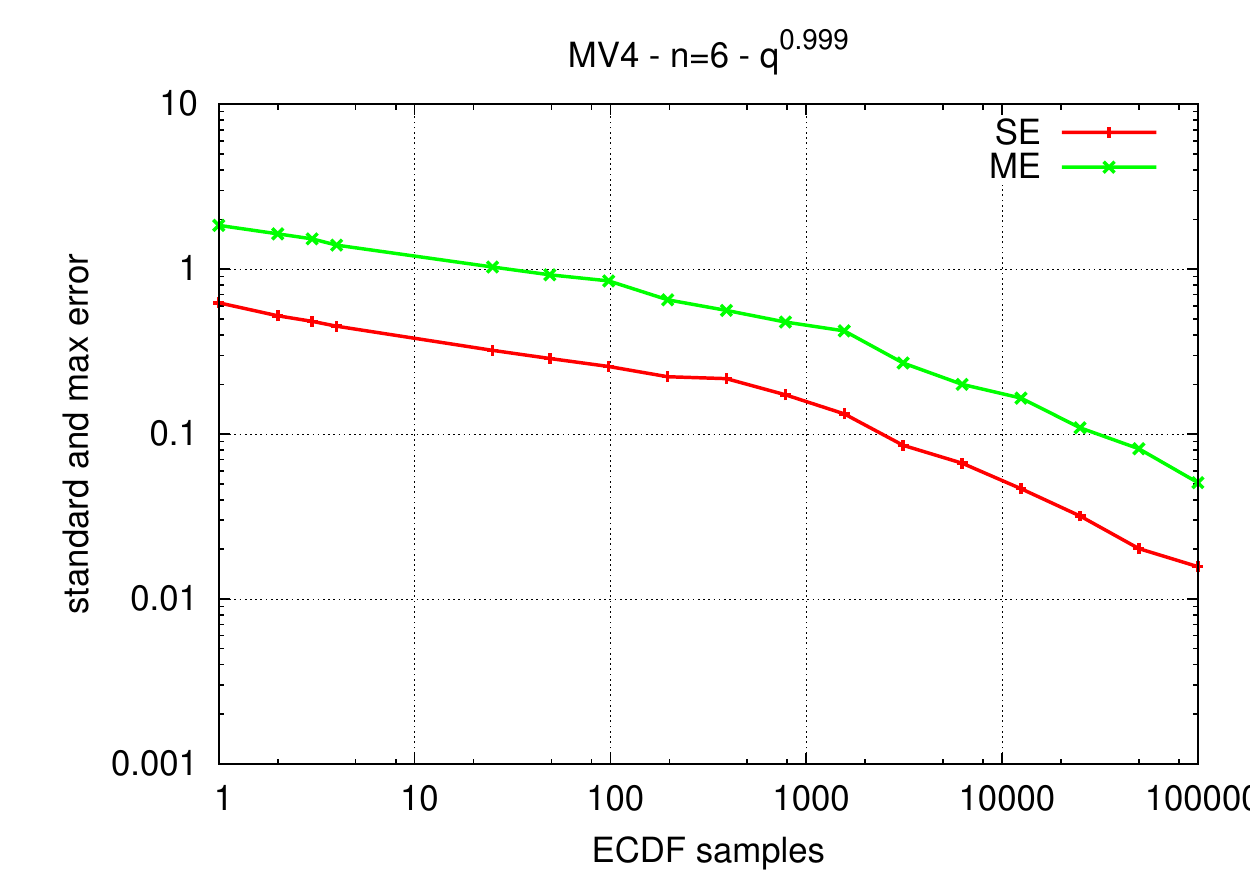}
}
\caption{Standard and maximum error for $q^{0.999}$ as a function of the number of samples for the MV4 ``MOST ROBUST'' solution with $n=6$.}
\label{error-vs-samples-ed-0999}
\end{figure}

\begin{figure}[h!]
\centerline{
\includegraphics[scale=0.8]{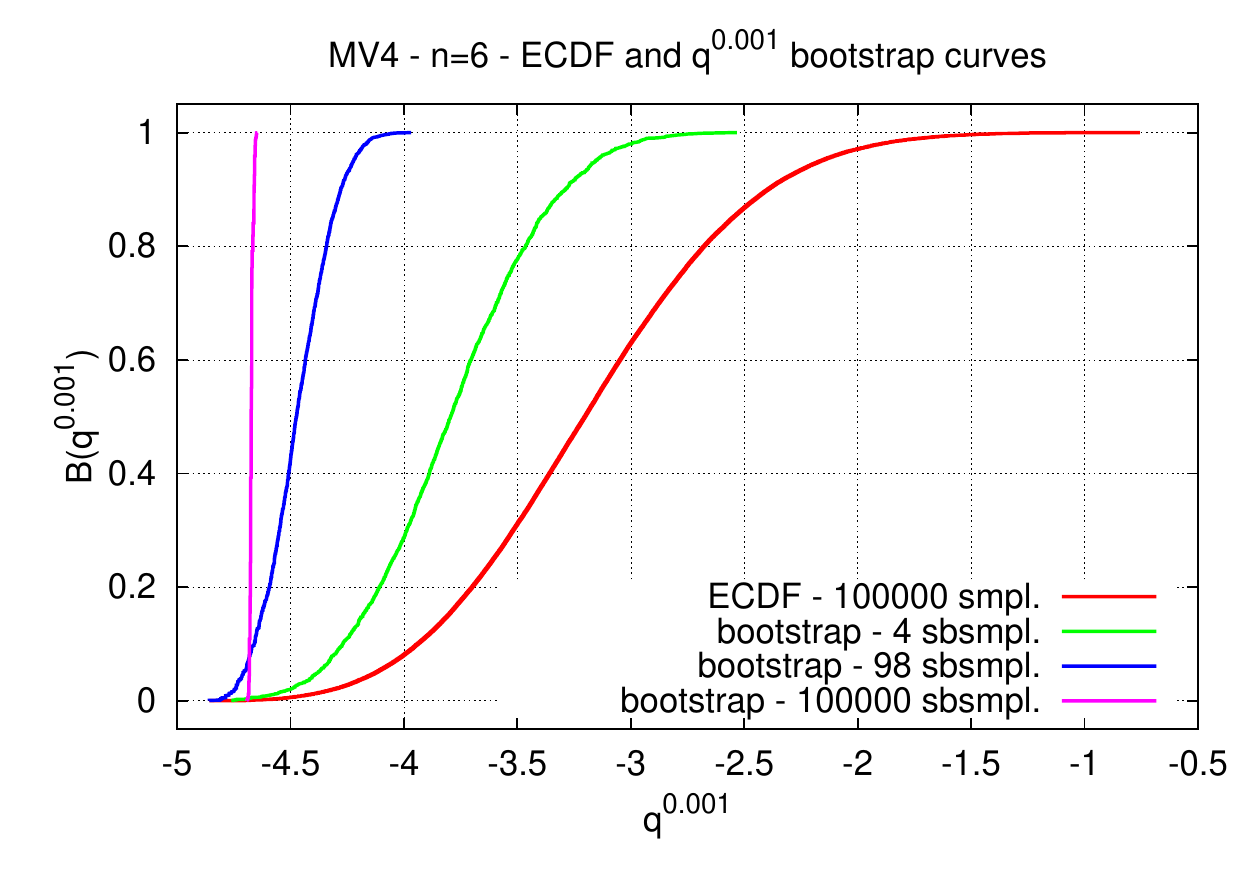}
}
\caption{Bootstrap curves for $q^{0.001}$ and high resolution ECDF for the MV4 ``MOST ROBUST'' solution with $n=6$.}
\label{ECDF+bootstrap-ed-0001}
\end{figure}

\begin{figure}[h!]
\centerline{
\includegraphics[scale=0.8]{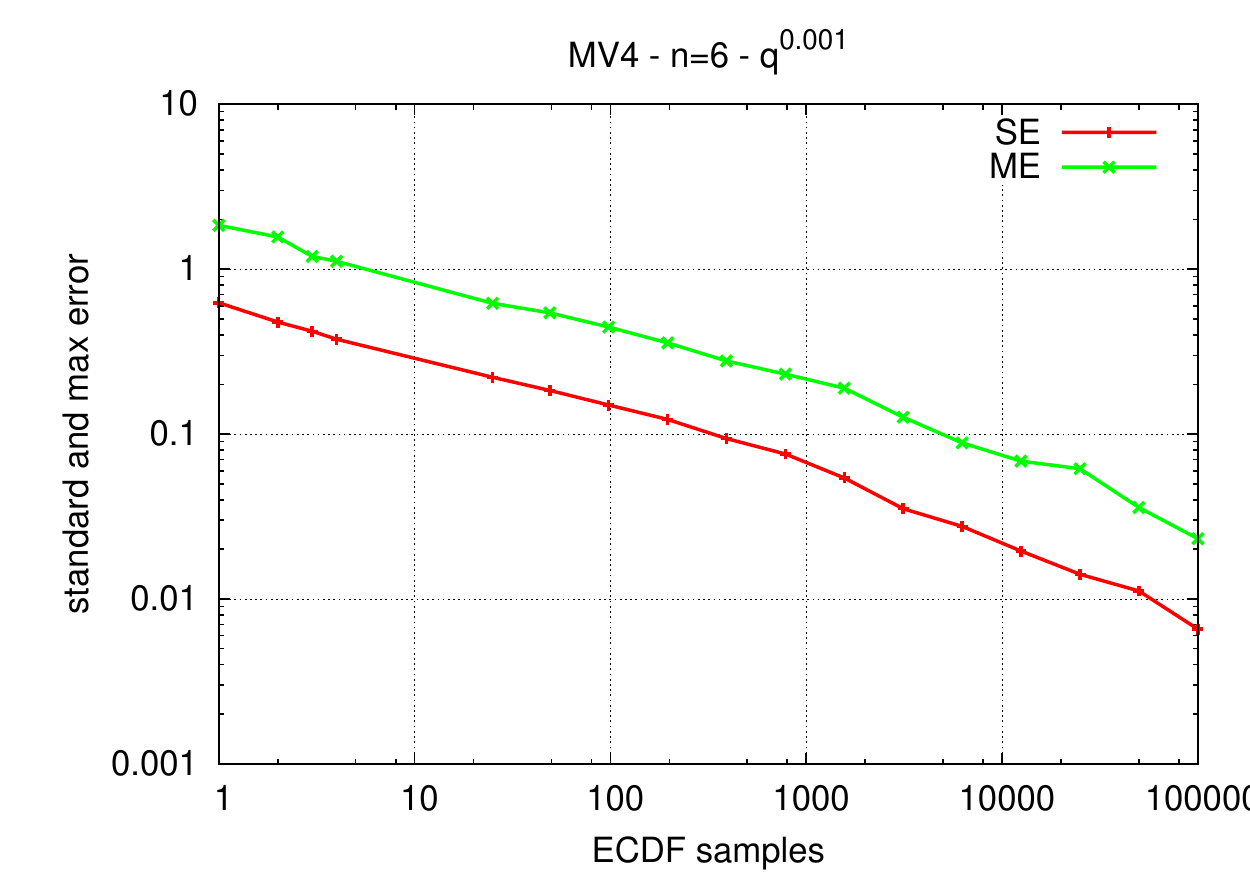}
}
\caption{Standard and maximum error for $q^{0.001}$ as a function of the number of samples for the MV4 ``MOST ROBUST'' solution with $n=6$.}
\label{error-vs-samples-ed-0001}
\end{figure}

The analysis of Figures
\ref{ECDF+bootstrap-ed-0999}--\ref{error-vs-samples-ed-0999},
related to $q^{0.999}$, and Figures
\ref{ECDF+bootstrap-ed-0001}--\ref{error-vs-samples-ed-0001}, related
to $q^{0.001}$, confirms the same trends previously detected, but
$\widehat{SE}$ and $\widehat{ME}$ errors remain higher, confirming the
fact that it is more difficult to describe the statistical tails compared
to median values.

\section{Many objectives formulation}

\begin{figure}[h!]
\centerline{
\includegraphics[scale=0.8]{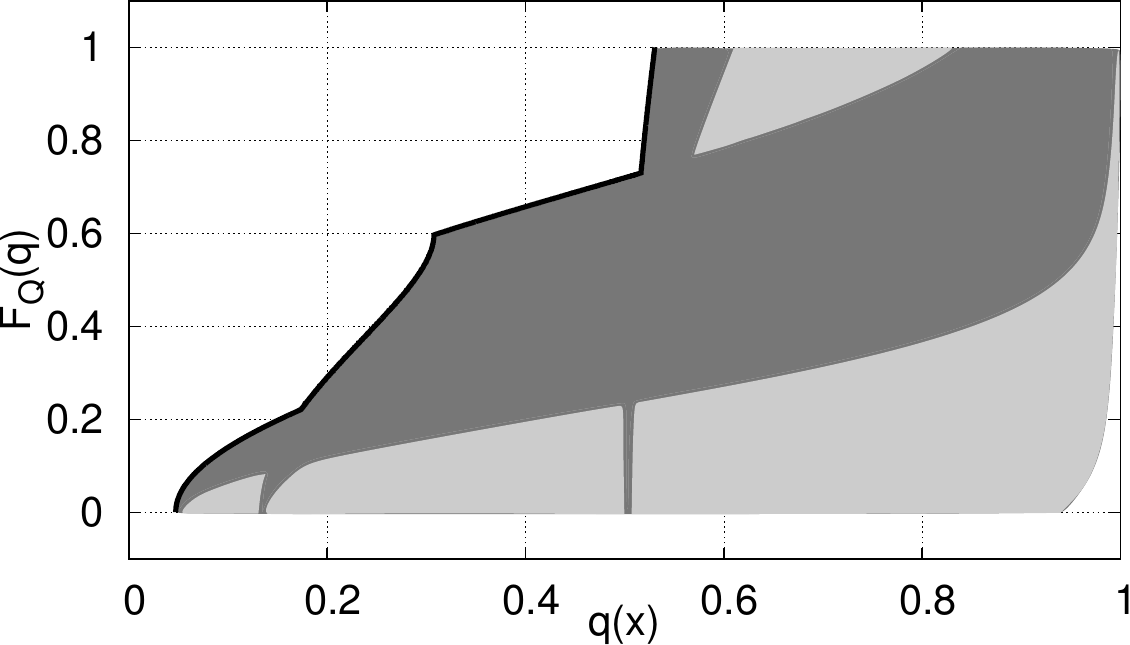}
}
\caption{Non-dominated envelop of all CDF curves.}
\label{envelop}
\end{figure}

The multi-objective formulation can be interpreted in a way that might
be advantageous in some cases. Suppose indeed to consider a very large
number, even infinite, of quantiles in the definition of the optimization
problem \ref{quantileopt}. This results in the problem that optimization
follows:
\begin{equation}
\left\{ {\begin{array}{*{20}c}
   {\mathop {\min}\limits_X q^{s_i }  = F_Q^{ - 1} \left( {s_i } \right)} \hfill & {} \hfill  \\
   {s_i  = i/n} \hfill & {i = 1,n} \hfill  \\
 \end{array} } \right.
\end{equation}
with $s \in [0,1]$ equally spaced probability thresholds. The aim of
this approach is to control and optimize the entire shape of the CDF.
This point can be illustrated with the help of an example. Consider again
the example described in Section \ref{illexam}. The CDF curves related to
all the possible elements of the set of definition of the function $q$,
numbered as (\ref{es}), are shown in light gray in Figure \ref{envelop}.
In the same figure, the CDF curves of the non dominated elements are
reported in dark gray and, finally, the envelope of the non dominated
CDF curves is reported in black.  It represents an ``IDEAL'' CDF for the
problem at hand, that is, a limit curve that outperforms any true CDF
and that can be adopted as a target in a template optimization process
as the one presented in \cite{Petrone2011b}.

\section{Solving an evidence theory optimization problem}

The same numerical machinery used for the GIDF can be immediately
adopted, with the addition of a few supplementary steps, to the
approximate estimation of the ``belief'' and ``plausibility'' curves
introduced in the ``Evidence Theory'' when this theory is applied
to solve robust optimization problems.  Indeed, belief, plausibility
and GIDF computations start from a sample of the uncertainty variable
parameter space, followed by the computation of the ECDF, to the heart
of which there is a sorting operation.  Similarly to the GIDF case,
the obtainable level of approximation for belief and plausibility curve
computation depends strongly on the quality and resolution of the ECDF
and is, hence, controllable.

The definition of the belief and plausibility is here biefly summarized
for the sake of clarity and completeness.  For each uncertain variable,
the evidence theory assigns, at given intervals, which can also be
partially overlapping, values called ``Basic Probability Assignment''
(BPA), which express the degree of confidence that the experts have on
that parameter. The BPA represents, in the framework of evidence theory,
the element that controls the propagation of information.  So, the BPA
allows to characterize each realization of a given random variable
$u$ in a way similar to the probability distribution function. The
conceptual difference is that here we are not dealing with true
probability distribution functions, but with probability bounds defined
by experts.  The BPA has to follow some simple rules that are reported
in the followings. Let $U$ be the set ov variation of an uncertain
parameter $u$, e.g. the interval $[-3, 5]$ in the example of Table
\ref{BPA}, and let $\Xi=\{E_1,\ldots,E_n\}$  a set of subsets of $U$
(e.g. $\bigcup_{i=1}^n E_i \subseteq U$).  A real, non-negative value
$m$, namely the basic probability, can be assigned to each element $E_i
\in \Xi$ according to the following rules:
\begin{equation}
\left\{
\begin{array}{ll}
m(E_i)\geq 0 & \forall E_i \in \Xi\\
m(\emptyset)=o\\
\displaystyle \sum_{E_i \in \Xi} m(E_i)=1
\end{array}
\right.
\end{equation}
The subset $\Phi=\{FE_1,\ldots,FE_m\} \subseteq \Xi$ of the elements with
non-zero value of $m$ is called the focal element set.  The extension
to more than one uncertain parameter is easily obtained with the help
of Cartesian product.  The assignment of a basic probability value
set to each uncertain parameter allows the computation of belief and
plausibility values for logical propositions involving these parameters.
Given a proposition $\alpha$, true in a subset $A$ of $U$, the belief
and plausibility values are defined as:
\begin{equation}\label{belandpl}
\begin{array}{rccl}
\mathrm{Bel}(\alpha)&=&\displaystyle\sum_{FE\subseteq A} & m(FE)\\
\mathrm{Pl} (\alpha)&=&\displaystyle\sum_{FE\cap A \neq \emptyset} & m(FE)
\end{array}
\end{equation}

\begin{table}
\begin{center}
\begin{tabular}{cc}
sub-interval & BPA \\ \hline
$[-5, -4]$ & 0.1\\
$[-3, 0]$ & 0.25\\
$[1, 3]$ & 0.65
\end{tabular}
\end{center}
\caption{BPA structure for an uncertain variable $u$ defined in $[-5, 3]$ interval.}
\label{BPA}
\end{table}

In the application of evidence theory to robust optimization, the
proposition $\alpha$ is often defined as a parametric inequality involving
a function $f(u,d)$ of uncertain parameters $u\in U$ and design variables
$d\in D$, and a threshold value $\nu$:
\begin{equation}
f(u,d)\leq \nu
\end{equation}
In these cases a strong connection is evident with the cumulative
distribution function, that will be used to obtain an approximate
evaluation of belief and plausibility curves.  The steps to follow are
enumerated in the following subsections

\subsubsection*{Affine mapping of the uncertain parameter space}

\begin{figure}
\centerline{
\includegraphics[scale=0.8]{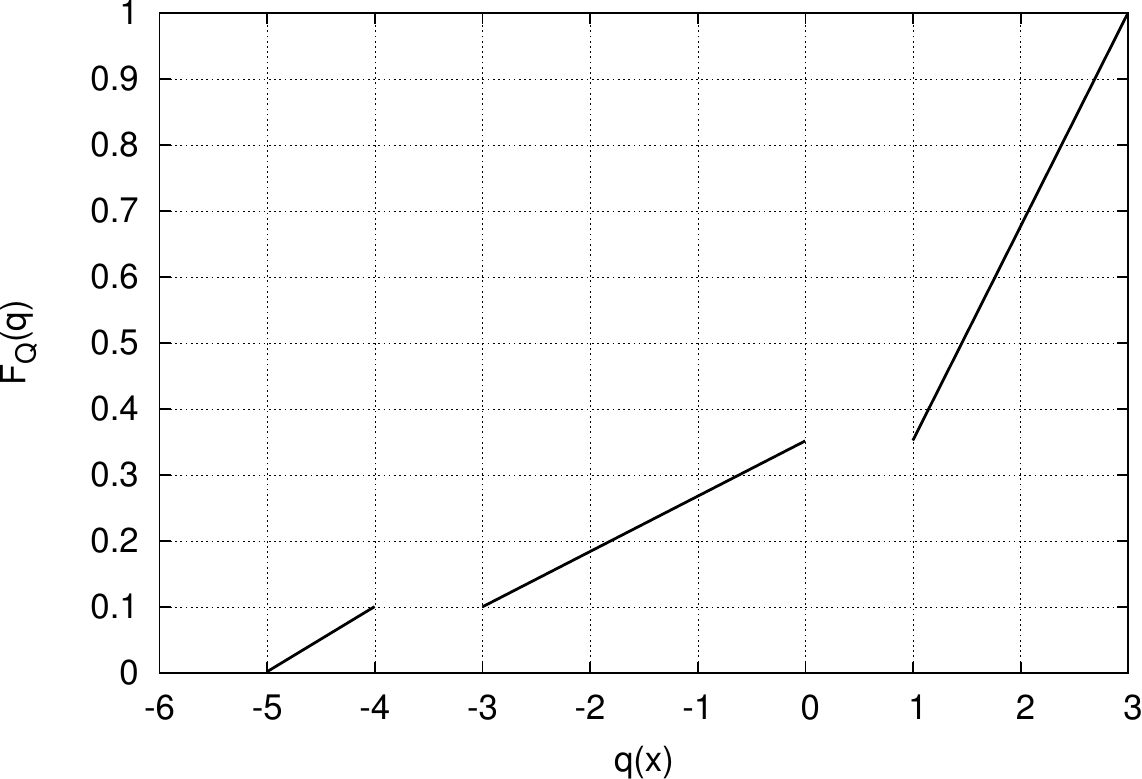}
}
\caption{Mapping of the BPA into a CDF.}
\label{BPAtoCDF}
\end{figure}

For each uncertain variable, the evidence theory assigns, at given
subsets, which can also be partially overlapping, the BPA values.
To the purpose of numerical computation, these BPA focus on assigned
regions of the design space that can be also non-connected. A sample BPA
definition for a given uncertain parameter is reported in table \ref{BPA}.
In order to restore the connectivity of the design space, and to make
easier the numerical computation, an affine mapping is applied to
the design space that is remapped in the unity cube.  In the present
method, the affine mapping is replaced by a CDF built for each uncertain
input variable.  This CDF assigns uniform selection probabilities to
each of the BPA intervals selected by the experts for that variable.
Figure \ref{BPAtoCDF} illustrates the mapping of the BPA defined in
Table \ref{BPA} into the CDF that replaces the affine mapping in the
current approach.

\subsubsection*{Calculation of belief and plausibility using the ECDF}

As it appears evident from formula \ref{belandpl}, the calculation of
$\mathrm{Bel}(f(u,d)\leq \nu)$, for a given $d$, requires, for each focal
element $FE \in A$, the computation of the maximum value of $f(u,d)$:
\begin{equation}
\max_{u\in FE}f(u,d)
\end{equation}
Similarly, the computation of the plausibility,
$\mathrm{Pl}(f(u,d)\leq \nu)$, requires the determination of
the set of minimum values of $f(u,d)$:
\begin{equation}
\min_{u\in FE}f(u,d)
\end{equation}

If a ECDF of $f(u,d)$ with $n$ samples is available, along with a mapping
that associates each point of the ECDF with its focal element $FE$, then
it is possible to estimate the belief curve. More formally, if $f_i$
is the $i$-th point of the ECDF, the mapping function can be written as:
\begin{equation}
p_i: f_i \in \mathrm{ECDF} \longmapsto \mathrm{FE}_i \in U
\end{equation}
and the related BPA value is indicated with $b_i$. {\bf This is the most
critical point.  The accuracy of this mapping is the key to obtain an
accurate estimation!}

For a given value $\nu^\star$ of the proposition threshold, the belief
value can be obtained considering the subset of the ECDF samples
defined by
\begin{equation}
M_{\nu^\star} = \left\{f_i: f_i>\nu^\star, f_i\in \mathrm{ECDF}\right\}
\end{equation}
Considering that the ECDF is, by definition, sorted with respect the $f_i$
values, it is convenient to define an index function $q=q(M_{\nu^\star})$
that maps the smaller element of $M_{\nu^\star}$ to the corresponding
element position in the ECDF.  The estimation of the belief is thus
given by:
\begin{equation}
\hat{B}_{\nu^\star}=\widehat{\mathrm{Bel}}(f(u,d)\leq \nu^\star)=1-\sum_{j=q}^{n}s_j
\end{equation}
where
\begin{equation}
s_j=\left\{
\begin{array}{cl}
b_j & \textrm{if} \; p_j\neq p_{\nu^\star} \; \textrm{and} \; p_j\neq p_k, \forall k=q,j-1\\
0   & \textrm{otherwise.}
\end{array}
\right.
\end{equation}
The calculation of plausibility is very similar, but with the important
difference that the BPA of the $FE_{\nu^\star}$ must be considered in
the related formulas.  Therefore let:
\begin{equation}
N_{\nu^\star} = \left\{f_i: f_i\geq\nu^\star, f_i\in \mathrm{ECDF}\right\}
\end{equation}
and $r=q(N_{\nu^\star})$.
So the plausibility estimation is given by
\begin{equation}
\hat{P}_{\nu^\star}=\widehat{\mathrm{Pl}}(f(u,d)\leq \nu^\star)=1-\sum_{j=r}^{n}t_j
\end{equation}
where
\begin{equation}
t_j=\left\{
\begin{array}{cl}
b_j & \textrm{if} \; p_j\neq p_k, \forall k=r,j-1\\
0   & \textrm{otherwise.}
\end{array}
\right.
\end{equation}

The approximated belief and plausibility curves thus obtained can be
immediately used in the optimization process.

\lstset{language=C, caption=C code for belief and plausibility estimation., label=BPEV, tabsize=4}
{\scriptsize
\begin{lstlisting}
#include <stdlib.h>
#include <stdio.h>
#include <math.h>
#include <string.h>

/*
 * Extended ECDF structure
 */
typedef struct {
	double f; /* Function value */
	char hash[100]; /* Focal point id label */
	int hash_min; /* In a vector of ecdf_t elements sorted with respect to f hash_min is set to 1
	                 for the first focal point with the same hash label encountered */
	int hash_max; /* In a vector of ecdf_t elements sorted with respect to f hash_max is set to 1
	                 for the last focal point with the same hash label encountered */
	double bpa; /* Basic probability assignment for the focal element to which f belongs */
	double belief;
	double plausibility;
} ecdf_t;

double MV1(double *d, double *u, int n);
double MV2(double *d, double *u, int n);
double MV8(double *d, double *u, int n);
void belief_and_plausibility_estimator(double (*f)(double*, double*, int),
		void (*affmap)(int, double*, double*), char (*feid)(double),
		double (*bpa_component)(char), double *d, double ul, double uu, int n,
		int m, ecdf_t *Ecdf);

void uranvec(double *x, int n, double xl, double xu) {
	int i;
	for (i = 0; i < n; i++) {
		x[i] = xl + (xu - xl) * drand48();
	}
}

/*
 * Begin of Problem Dependent Data
 */
double affine_mapping_MV1_comp(double x) {
	/* Affine mapping for the MV1 function variable vector components */
	double y;
	if (x >= 0 && x <= 0.1)
		y = -5 + (-4 + 5) * x / 0.1;
	else if (x > 0.1 && x <= 0.1 + 0.25)
		y = -3 + (0 + 3) * (x - 0.1) / 0.25;
	else
		y = 1 + (3 - 1) * (x - 0.35) / 0.65;
	return y;
}

void affine_mapping_MV1(int n, double *x, double *mx) {
	int i;
	for (i = 0; i < n; i++)
		mx[i] = affine_mapping_MV1_comp(x[i]);
}

char focal_element_id_MV1(double x) {
	/*
	 * Computes the focal element id related to a single vector variable component
	 * (after the affine transformation) for the MV1 function.
	 */
	char y;
	if (x >= -5 && x <= -4)
		y = '0';
	else if (x > -3 && x <= 0)
		y = '1';
	else if (x > 1 && x <= 3)
		y = '2';
	else
		exit(0); /* ERROR: Out of considered intervals. */
	return y;
}

double bpa_component_MV1(char x) {
	/*
	 * Basic probability assignment for a single vector variable component
	 * of the MV1 function.
	 */
	int iu;
	double bpa_c = 0;
	iu = x - '0';
	switch (iu) {
	case 0:
		bpa_c = 0.10;
		break;
	case 1:
		bpa_c = 0.25;
		break;
	case 2:
		bpa_c = 0.65;
		break;
	}
	return bpa_c;
}

double MV1(double *d, double *u, int n) {
	/*
	 * $\sum_{i=1}^{n} d_i u_i^2, d\in[1,5]^n, u\in[-5,3]^n
	 */
	double objective;
	int i, j;
	objective = 0;
	for (j = 1; j <= n; j++) {
		i = j - 1;
		objective += d[i] * (u[i] * u[i]);
	}
	return objective;
}

double MV2(double *d, double *u, int n) {
	/*
	 * $\sum_{i=1}^{n}(u_i-d_i)^2, d\in[1,5]^n, u\in[-5,3]^n
	 */
	double objective;
	int i, j;
	objective = 0;
	for (j = 1; j <= n; j++) {
		i = j - 1;
		objective += (u[i] - d[i]) * (u[i] - d[i]);
	}
	return objective;
}

double MV8(double *d, double *u, int n) {
	/*
	 * $\sum_{i=1}^{n}(2\pi-u_i)\cos(u_i-d_i)^2, d\in[0,3]^n, u\in[0,2\pi]^n
	 */
	double objective;
	int i, j;
	objective = 0;
	for (j = 1; j <= n; j++) {
		i = j - 1;
		objective += (2 * M_PI - u[i]) * cos(u[i] - d[i]);
	}
	return objective / n;
}
/*
 * End of Problem Dependent Data
 */

static int dcomp_ecdf_t(ecdf_t *a, ecdf_t *b) {
	if (a->f < b->f)
		return -1;
	if (a->f > b->f)
		return +1;
	return 0;
}

void belief_and_plausibility_estimator(double (*f)(double*, double*, int),
		void (*affmap)(int, double*, double*), char (*feid)(double),
		double (*bpa_component)(char), double *d, double ul, double uu, int n,
		int m, ecdf_t *Ecdf) {
	double u[n], amu[n];
	int i, j;
	for (i = 0; i < m; i++) {
		uranvec(u, n, ul, uu);
		affmap(n, u, amu);
		Ecdf[i].f = f(d, amu, n);
		Ecdf[i].bpa = 1;
		for (j = 0; j < n; j++) {
			Ecdf[i].hash[j] = feid(amu[j]);
			Ecdf[i].bpa *= bpa_component(Ecdf[i].hash[j]);
		}
		Ecdf[i].hash[n] = 0;
		Ecdf[i].hash_min = 1;
		Ecdf[i].hash_max = 1;
	}
	qsort(Ecdf, m, sizeof(ecdf_t),
			(int (*)(const void *, const void *)) dcomp_ecdf_t);
	for (i = m - 1; i >= 0; i--) {
		if (Ecdf[i].hash_max == 1)
			for (j = i - 1; j >= 0; j--) {
				if (strcmp(Ecdf[i].hash, Ecdf[j].hash) == 0)
					Ecdf[j].hash_max = 0;
			}
	}
	for (i = 0; i < m; i++) {
		if (Ecdf[i].hash_min == 1)
			for (j = i + 1; j < m; j++) {
				if (strcmp(Ecdf[i].hash, Ecdf[j].hash) == 0)
					Ecdf[j].hash_min = 0;
			}
	}
	for (i = 0; i < m; i++) {
		Ecdf[i].belief = 1;
		for (j = i + 1; j < m; j++) {
			if (Ecdf[j].hash_max == 1)
				Ecdf[i].belief -= Ecdf[j].bpa;
		}
	}

	for (i = 0; i < m; i++) {
		Ecdf[i].plausibility = 1;
		for (j = i + 1; j < m; j++) {
			if (Ecdf[j].hash_min == 1)
				Ecdf[i].plausibility -= Ecdf[j].bpa;
		}
	}
}

#define NUMBER_OF_OUTPUT_VARIABLES 10000
#define VECTOR_VARIABLE_SIZE 1

int main() {
	/*
	 * Compute belief and plausibility
	 */
	ecdf_t Ecdf[NUMBER_OF_OUTPUT_VARIABLES];
	double r[VECTOR_VARIABLE_SIZE];
	int n = VECTOR_VARIABLE_SIZE, i, m = NUMBER_OF_OUTPUT_VARIABLES;
	FILE *BEL, *PLA;
	srand48(1234);
	for (i = 0; i < n; i++)
		r[i] = 1;
	belief_and_plausibility_estimator(MV1, affine_mapping_MV1,
			focal_element_id_MV1, bpa_component_MV1, r, 0, 1, n,
			NUMBER_OF_OUTPUT_VARIABLES, Ecdf);
	BEL = fopen("belief.dat", "w");
	for (i = 0; i < m; i++) {
		fprintf(BEL, "HASH %s %d %d %lf %lf %lf\n", Ecdf[i].hash,
				Ecdf[i].hash_min, Ecdf[i].hash_max, Ecdf[i].f, Ecdf[i].bpa,
				Ecdf[i].belief);
	}
	fclose(BEL);
	PLA = fopen("plausibility.dat", "w");
	for (i = 0; i < m; i++) {
		fprintf(PLA, "HASH %s %d %d %lf %lf %lf\n", Ecdf[i].hash,
				Ecdf[i].hash_min, Ecdf[i].hash_max, Ecdf[i].f, Ecdf[i].bpa,
				Ecdf[i].plausibility);
	}
	fclose(PLA);
	return 0;
}

\end{lstlisting}
}

\subsection{A simple example}

The MV1 test function (see ref. \cite{Vasile2011}) is used to illustrate
the ECDF based belief and plausibility estimation method:
\begin{equation}
f=\sum_{i=1}^n d_i u_i^2
\end{equation}
with design parameter vector $\mathbf{d}\in[1,5]^n$ and uncertain
parameter vector $\mathbf{u}\in[-5,3]^n$.  The BPA is the one already
reported in Table \ref{BPA}, and the computations have been made with
the design parameter vector assigned to $\mathbf{d}=[1]^n$.

\begin{figure}
\centerline{
\includegraphics[scale=1.0]{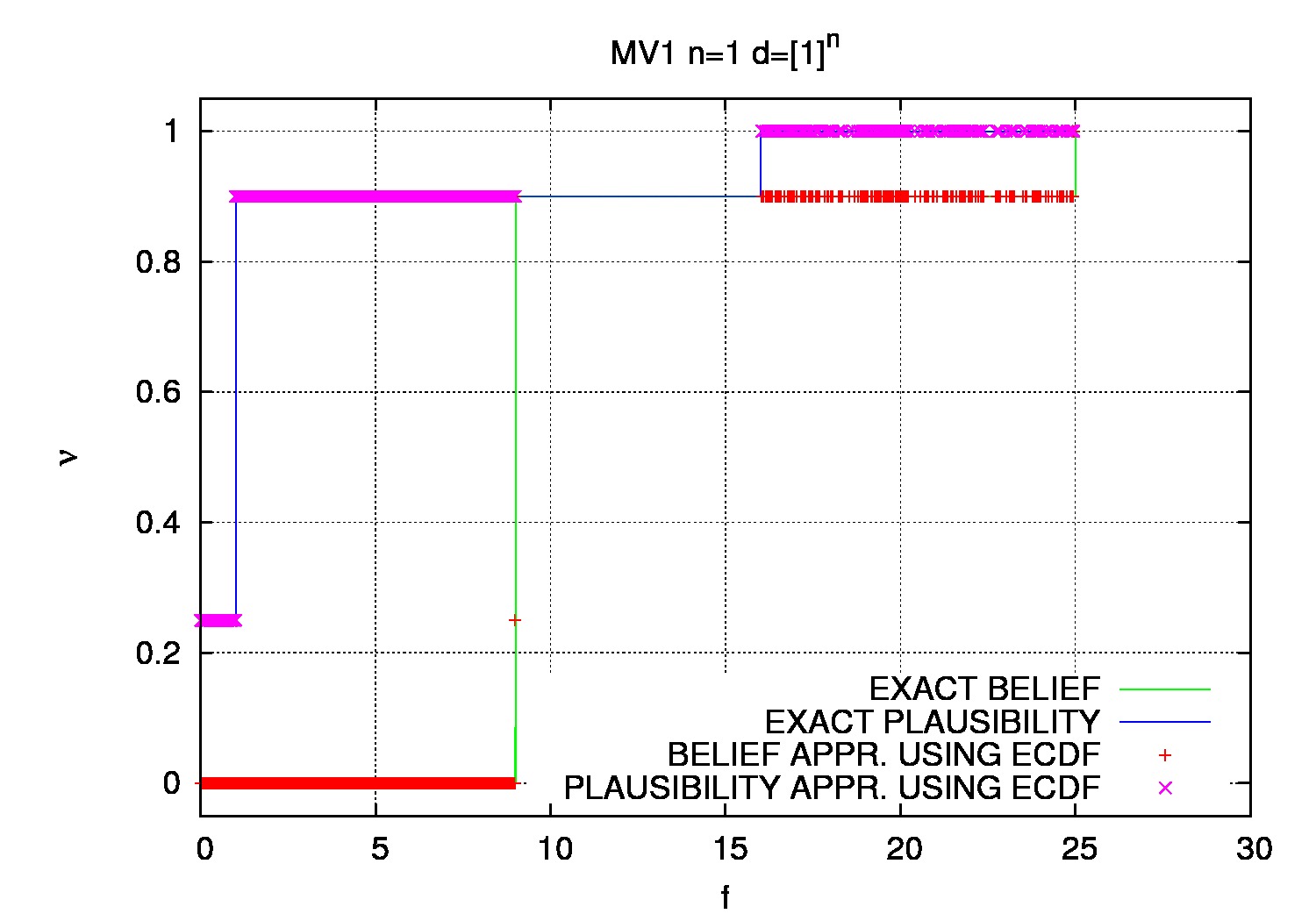}
}
\caption{MV1 test function with design and uncertain parameter space size equal to 1.}
\label{MV1}
\end{figure}

\begin{figure}
\centerline{
\includegraphics[scale=1.0]{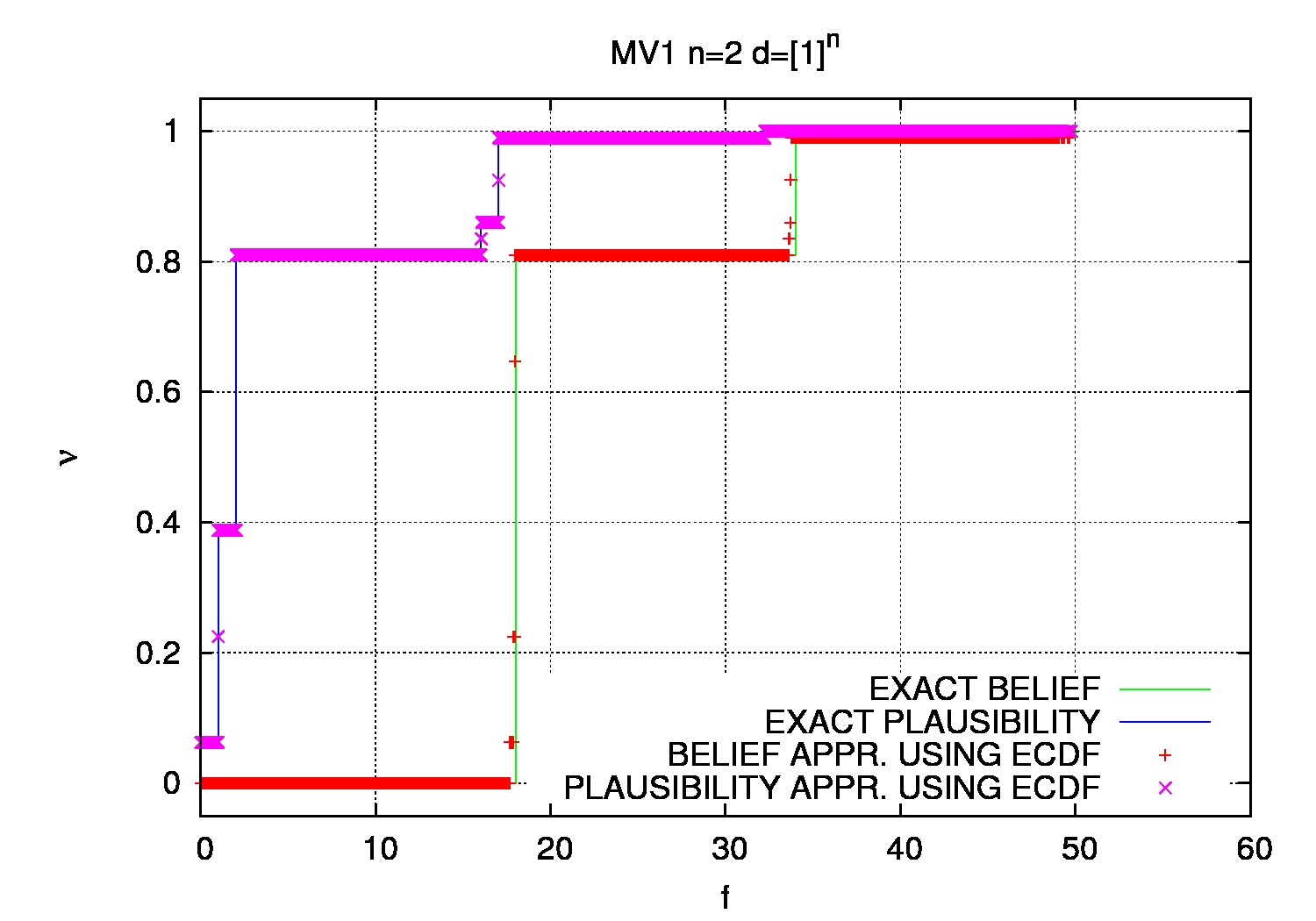}
}
\caption{MV1 test function with design and uncertain parameter space size equal to 2.}
\label{MV1-2}
\end{figure}

\begin{figure}
\centerline{
\includegraphics[scale=1.0]{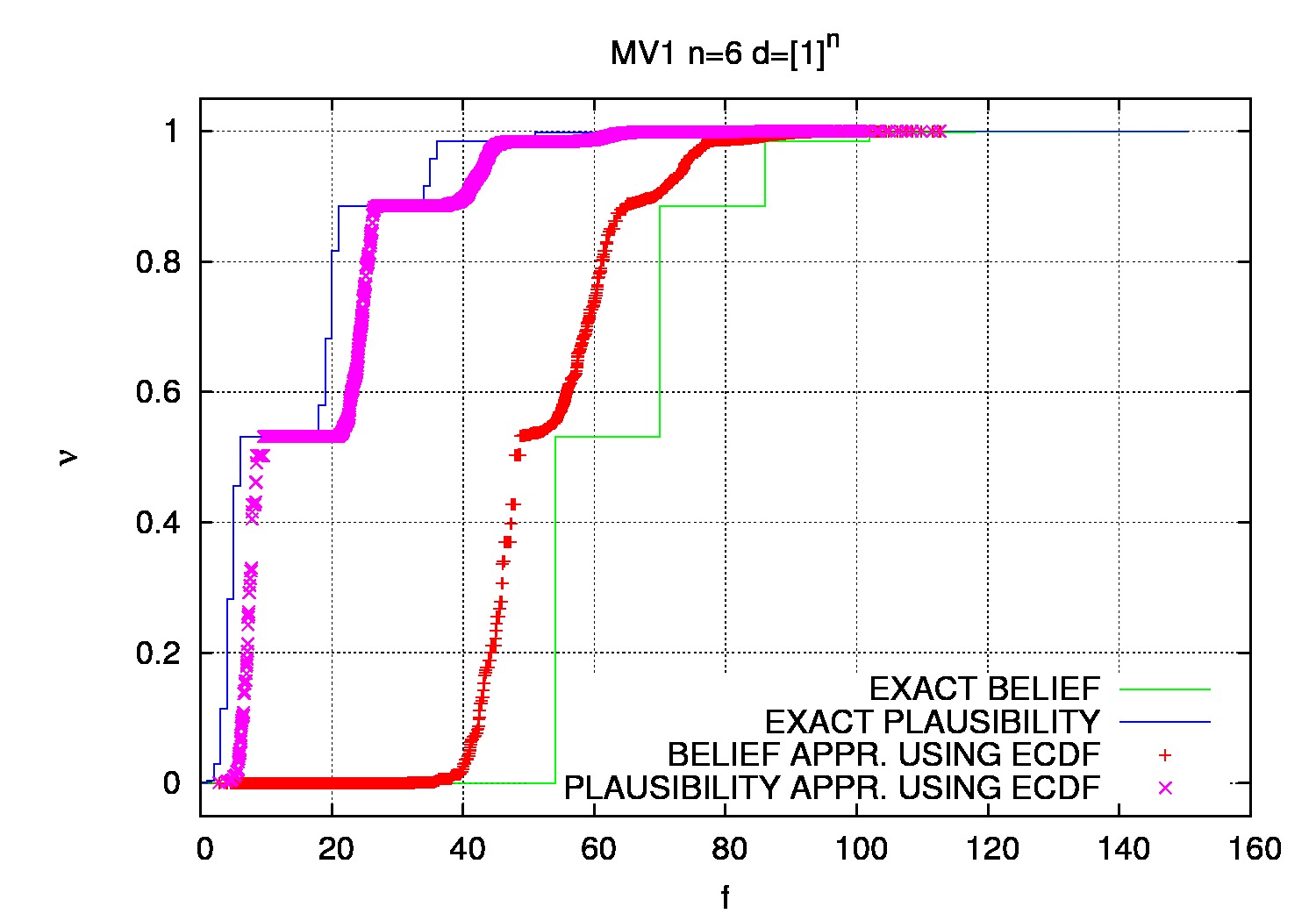}
}
\caption{MV1 test function with design and uncertain parameter space size equal to 6.}
\label{MV1-6}
\end{figure}

The curves of belief and plausibility for the test function MV1,
calculated using the method above, are shown in Figures \ref{MV1},
\ref{MV1-2} and \ref{MV1-6}. These figures show also the exact belief
and plausibility curves for the given test case in order to allow the
evaluation of the actual potential of the method. What is clear is
that the obtained estimate, while acceptable, is not conservative,
and therefore some caution is required in the use of the results.

\section{Conclusions}

An alternative approach to the optimization under uncertainty has been
introduced with the help some simple test cases. The first one was a
very simple and well behaved test function that had the only purpose
of illustrating the main feature of the quantile based optimization
method. A further example function with variable number of uncertain
and design parameters was then introduced to highlight the critical
points of the method related to the problem dimension. In particular,
it has been discussed and illustrated how the problem features change
when the number of random variables involved increases.  Work is under
way to better characterize the error introduced by the ECDF estimator
as a function of the number of random parameters.

The algorithm used for optimization is a classical genetic algorithm,
but, to further improve the proposed procedure, an optimization algorithm
capable of accounting for the errors in the estimation of the CDF has
to be conceived. This is a very important topic and it will be subject
of next research work. In particular, to reduce the curse of dimensionality the effect of different sampling methodologies, like stochastic collocation, on the estimation of the ECDF will be considered in future works. Indeed the possibility to use the error on the ECDF estimator to properly refine the probability space using adaptive uncertainty quantification algorithms will be explored.

It is, finally, very important to to relate the results of this new
optimization approach to those deriving from the application of more
conventional methods, and to introduce a quantitative approach when
different algorithms for robust optimization are compared.

\section{Appendix}

\bibliographystyle{spmpsci}
\bibliography{ro-biblio,bibliogr}

\begin{thebibliography}{10}
\providecommand{\url}[1]{{#1}}
\providecommand{\urlprefix}{URL }
\expandafter\ifx\csname urlstyle\endcsname\relax
  \providecommand{\doi}[1]{DOI~\discretionary{}{}{}#1}\else
  \providecommand{\doi}{DOI~\discretionary{}{}{}\begingroup
  \urlstyle{rm}\Url}\fi

\bibitem{Bratley1988}
Bratley, P., Fox, B.L.: Algorithm 659: Implementing sobol's quasirandom
  sequence generator.
\newblock ACM Trans. Math. Softw. \textbf{14}(1), 88--100 (1988).
\newblock \doi{10.1145/42288.214372}.
\newblock \urlprefix\url{http://doi.acm.org/10.1145/42288.214372}

\bibitem{Deb:2001}
Deb, K.: Multi-Objective Optimization Using Evolutionary Algorithms.
\newblock John Wiley \& Sons, Inc., New York, NY, USA (2001).
\newblock ISBN:047187339X

\bibitem{Diaconis1983}
Diaconis, P., Efron, B.: Computer intensive methods in statistics.
\newblock Tech. Rep.~83, Division Of Biostatistics, Stanford University (1983)

\bibitem{Efron1979}
Efron, B.: Bootstrap methods: Another look at the jackknife.
\newblock Annals of Statistics \textbf{7}(1), 1--26 (1979).
\newblock \urlprefix\url{http://projecteuclid.org/euclid.aos/1176344552}

\bibitem{coelho2011multi}
Filomeno~Coelho, R., Bouillard, P.: Multi-objective reliability-based
  optimization with stochastic metamodels.
\newblock Evolutionary Computation \textbf{19}(4), 525--560 (2011)

\bibitem{hughes2001evolutionary}
Hughes, E.: Evolutionary multi-objective ranking with uncertainty and noise.
\newblock In: Evolutionary multi-criterion optimization, pp. 329--343. Springer
  (2001)

\bibitem{Neumann1953}
von Neumann, J., Morgenstern, O.: Theory Of Games And Economic Behavior.
\newblock Princeton University Press, Princeton (1953)

\bibitem{Petrone2011b}
Petrone, G., Iaccarino, G., Quagliarella, D.: Robustness criteria in
  optimization under uncertainty.
\newblock In: C.~Poloni, D.~Quagliarella, J.~Periaux, N.~Gauger,
  K.~Giannakoglou (eds.) EUROGEN 2011 PROCEEDINGS --- Evolutionary and
  Deterministic Methods for Design, Optimization and Control with Applications
  to Industrial and Societal Problems, ECCOMAS Thematic Conference, pp.
  244--252. ECCOMAS, CIRA, Capua, Italy (2011)

\bibitem{poloni2004robust}
Poloni, C., Padovan, Parussini, L., Pieri, S., Pediroda, V.: Robust design of
  aircraft components: a multi-objective optimization problem.
\newblock In: Von Karman Institute for Fluid Dynamics, Lecture Series (2004).
\newblock Von Karman Institute for Fluid Dynamics, Lecture Series 2004-07

\bibitem{Quaglia:2000b}
Quagliarella, D., Vicini, A.: {GAs} for aerodynamic shape design {II}:
  multiobjective optimization and multi-criteria design.
\newblock In: Lecture Series 2000-07, Genetic Algorithms for Optimisation in
  Aeronautics and Turbomachinery. Von Karman Institute, Belgium (2000)

\bibitem{Serfling2008}
Serfling, R.J.: Approximation Theorems of Mathematical Statistics.
\newblock John Wiley \& Sons, Inc. (2008).
\newblock \doi{10.1002/9780470316481}.
\newblock \urlprefix\url{http://dx.doi.org/10.1002/9780470316481.indsub}

\bibitem{Sobol1994}
Sobol, I.M.: A Primer for the Monte Carlo Method.
\newblock CRC Press (1994)

\bibitem{teich2001pareto}
Teich, J.: Pareto-front exploration with uncertain objectives.
\newblock In: Evolutionary multi-criterion optimization, pp. 314--328. Springer
  (2001)

\bibitem{Vaart1998}
van~der Vaart, A.W.: Asymptotic Statistics.
\newblock Cambridge University Press (1998).
\newblock \urlprefix\url{http://dx.doi.org/10.1017/CBO9780511802256}

\bibitem{Vasile2011}
Vasile, M., Minisci, E.: An evolutionary approach to evidence-based
  multi-disciplinary robust design optimisation.
\newblock In: C.~Poloni, D.~Quagliarella, J.~Periaux, N.~Gauger,
  K.~Giannakoglou (eds.) EUROGEN 2011 PROCEEDINGS --- Evolutionary and
  Deterministic Methods for Design, Optimization and Control with Applications
  to Industrial and Societal Problems, ECCOMAS Thematic Conference. ECCOMAS,
  CIRA, Capua, Italy (2011)

\bibitem{Vicini:97b}
Vicini, A., Quagliarella, D.: Inverse and direct airfoil design using a
  multiobjective genetic algorithm.
\newblock AIAA Journal \textbf{35}(9), 1499--1505 (1997)

\bibitem{CDF}
Wikipedia: Cumulative distribution function.
\newblock
  \urlprefix\url{http://en.wikipedia.org/wiki/Cumulative\_distribution\_function}

\bibitem{ECDF}
Woo, C.: Empirical distribution function (version 4).
\newblock
  \urlprefix\url{http://planetmath.org/EmpiricalDistributionFunction.html}

\end{thebibliography}

\end{document}